\documentclass[12pt,a4paper]{article}

\usepackage[cp1251]{inputenc}   % Windows
\usepackage[english]{babel}

\newtheorem{teoen}{Theorem}
\newtheorem{lemen}{Lemma}

\usepackage{amsmath}
\usepackage{amsfonts,amssymb}

\usepackage{graphicx}

\usepackage{psfrag}                    
\usepackage{colordvi}

\usepackage{cite}

\usepackage{srcltx}

\makeatletter
\def\thmstyle{\it} % style of text in theorem environment
\def\@begintheorem#1#2{\it \trivlist \item[\hskip
        \labelsep{\bf #1\ #2.}]\thmstyle}
\def\@opargbegintheorem#1#2#3{\it \trivlist \item[\hskip
        \labelsep{\bf #1\ #2\ (#3).}]\thmstyle}
\makeatother

\begin{document}

\begin{center}
{\bf
The Upper Bound for the Lebesgue Constant for Lagrange Interpolation
in Equally Spaced Points of the Triangle
\\[2ex]
N.\,Baidakova}\\[3ex]
\end{center}

\begin{center}
\small{\it School of Mathematical Science, Tel-Aviv University, Israel}\\
%\small{\it School of Mathematical Science, Tel-Aviv University, Ramat Aviv, Tel Aviv, 69978, Israel}
\small{\it nataliabdkv@gmail.com}
\end{center}

{\small
An upper bound for the Lebesgue constant (the supremum norm)
of the operator of interpolation of a function
in equally spaced points of a triangle by a polynomial of total degree
less than or equal to $n$
is obtained.
%, preserving a known order in $n$ as $n\to\infty$.
Earlier, the rate of increase of the Lebesgue constants
with respect to $n$ for an arbitrary $d$-dimensional
simplex was established by the author. The explicit
upper bound proved in this article refines this result
for $d=2$.
%, preserving the known order.
}

%{\small
%An upper bound
%for the Lebesgue constant (the supremum norm)
%of the operator for interpolating a function
%over equidistant nodes of a triangle by a polynomial of degree $n$
%is obtained, preserving a known order in $n$ as $n\to\infty$.
%}

\vspace{5mm}

{\small

}

\vspace{5mm}

 \bigskip

\section{Introduction}

Let $\Delta$ be a non-degenerate simplex in $\mathbb{R}^d$
with vertices $\bar{a}_1,\ldots,\bar{a}_{d+1}$.
Any point $u\in\Delta$ can be written as
$$
u=\sum\limits_{r=1}^{d+1}
\lambda_r\bar{a}_r,
\quad
\sum\limits_{r=1}^{d+1}
\lambda_r=1,
\quad
\lambda_r\ge 0
\quad\mbox{for}\quad
1\le r\le d+1.
$$
The elements of
the tuple $\lambda=(\lambda_1,\ldots,\lambda_{d+1})$
are called the barycentric coordinates of the point $u$.
Further we will identify $u$ with
$
\lambda
\overset{def}{=}
(\lambda_1,\ldots,\lambda_{d+1})
\in \Delta.
$
Let's denote
$$
I=\left\{
i=(i_1,\ldots,i_{d+1}):\quad
i_1,\ldots,i_{d+1}\in \mathbb{Z}_+,
\quad
\sum\limits_{r=1}^{d+1}i_r=n
\right\},
$$
where $\mathbb{Z}_+$ is the set of non-negative integers.
Let $P_n^d=P_n^d[f]=P_n^d[f](u) $ be a polynomial
of total degree
less than or equal to
$n
%\in \mathbb{N}
$ interpolating
the values of the function $f\in C( \Delta)$ at
equidistant nodes $a_i$ of the simplex $\Delta$ , where
$$
a_i
=\left(\frac{i_1}{n},\ldots,\frac{i_{d+1}}{n}\right),
\quad
i=(i_1,\ldots,i_{d+1})\in I.
$$
If we take the polynomials $l_{j}$
of total degree
less than or equal to  $n$ defined by the conditions
$$
l_j(a_i)
=\delta_{i,j}
=
\left\{
\begin{array}{l}
1,\quad i=j,\\
0,\quad i\neq j,
\end{array}
\right.
\quad i,j\in I,
$$
then
$$
P_n^d(u)
=
\sum\limits_{i\in I}f(a_i)l_i(u).
$$
The polynomials $l_{j}$
are called Lagrange's fundamental interpolating polynomials.

Denote by
${\mathcal{L}}_n^d={\mathcal{L}}_n^d(u)={\mathcal{L}}_n^d(\lambda)$
the Lebesgue function of the Lagrange interpolation process
at equally spaced nodes of the simplex, i.e., the norm of a functional
in the space  $C(\Delta)$ that assigns
to each continuous function $f\in C(\Delta)$ the value of its
interpolation polynomial at the point $u\in\Delta$.
By $\Lambda_n^d$ we denote the Lebesgue constant
of this interpolation process,
i.e., the norm of an operator from $C(\Delta)$ to  $C(\Delta)$
which assigns to each continuous function its interpolation polynomial $P_n^d[f]$.
Thus
$$
{\mathcal{L}}_n^d(u)=
{\mathcal{L}}_n^d(\lambda)=\sup\limits_{f\in C(\Delta) \atop {f\neq 0}}
\frac{|P_n^d(u)|}{\|f\|_{C(\Delta)}}~,
\quad
\quad
{\displaystyle
\Lambda_n^d=\sup\limits_{f\in C(\Delta)\atop {f\neq 0}}
\frac{\|P_n^d\|_{C(\Delta)}}{\|f\|_{C(\Delta)}}}~.
$$
It is known that
\begin{equation}
\label{const_Leb}
\Lambda_n^d
=\max\limits_{u\in\Delta}{\mathcal{L}}_n^d(u)
=\max\limits_{\lambda\in\Delta}{\mathcal{L}}_n^d(\lambda),
\quad
\end{equation}
\begin{equation}
\label{funcLeb}
{\mathcal{L}}_n^d(\lambda)
=
\sum\limits_{i\in I}|l_i(\lambda)|.
\end{equation}

For $d=1$, in 1940 Turetskii~\cite{Ture1940} proved
the following asymptotics of the Lebesgue constant:
\begin{equation}
\label{Turec_vved}
\Lambda_n^1=\frac{2^{n+1}}{e~n\ln n}(1+\varepsilon_n), \quad
\mbox{where}\quad \lim\limits_{n\to \infty}\varepsilon_n=0.
\end{equation}
This result was also published in his book of 1968~\cite{Ture1968}.
In~1961,
%equality~
(\ref{Turec_vved}) was independently
proved again by Sch\"{o}nhage~\cite{Scho11961}
(more precisely, the formulation of the result in~\cite{Scho11961}
is similar to the formula~(\ref{Turec_vved}); however, 
a slightly stronger result is proved in~\cite{Scho11961},
where the first terms of the asymptotic expansion
of the Lebesgue constant~$\Lambda_n ^1$ for $n\to\infty$ were given).
A review of other results
%obtained up to 1991~
can be found in the article by Trefethen and Weideman~\cite{Tref1991}.
In~1992, Mills and Smith~\cite{MillSmit1992} obtained an asymptotic
expansion for $\ln \Lambda_n^1$, which refines the result from~\cite{Scho11961}.
In 2004, Eisinberg, Fedele, Franz~\cite{EisiFedeFran2004}
got a new asymptotic formula for~$\Lambda_n^1$
and numerically demonstrated the advantage of the formula
for $33\le n\le 200$ in terms of the relative error.

In 1983, for arbitrary nonnegative integers $d,n$
Bos~ \cite{Bos1983} proved that
$$
\Lambda_n^d\le
\left(
\begin{array}{c}
2n-1\\
n
\end{array}
\right)
\asymp\frac{4^n}{\sqrt{n}}~;
\quad
\Lambda_n^d\to
\left(
\begin{array}{c}
2n-1\\
n
\end{array}
\right)
\quad
\mbox{if}
\quad
d\to \infty.
$$
In 1988, Bloom~\cite{Bloo1988} established that
$$
\lim\limits_{n\to\infty}
\frac{\ln \Lambda_n^d}{n}
=
\ln 2,
$$
and in 2005, the author~\cite{Baid2005} proved that
\begin{equation}
\label{order}
C_1
\frac{2^n}{n\ln n}
\leq
\Lambda_n^d
\leq
C_2(d)
\frac{2^n}{n\ln n}~,
\end{equation}
where $C_1$ is some positive constant,
and $C_2(d)$ may depend on~$d$ (but does not depend on $n$).

The purpose of this paper is to obtain an upper bound
of $\Lambda_n^d$ for $d=2$
with preservation of the known rate of increase with respect to $n$, i.e. essentially
to estimate the constant $C_2(d)$ in~(\ref{order}) for $d=2$.
Namely, for $n\ge 4$, we prove that
$$
\Lambda_n^2
\le
(7+\mu_n)
\frac{2^{n+1}}{en(\ln n -\ln 2)}
\left(1+\frac{15}{n-3}\right),
\quad
\mbox{where}
\quad
\mu_n
\le
\frac{4en(\ln n)^3}{2^n/3}
+
\frac{en^2\ln n}{2^n}
\to 0\quad
\mbox{as}
\quad
n\to \infty.
$$

The paper is organized as follows.
In Sections $2-4$, some auxiliary results are given.
In Sections $5-7$, the problem of estimating $\Lambda_n^2$
is reduced to the problem of estimating the function ${\mathcal{L}}_n^2(\lambda)$
for the case $0\le \lambda_2,\lambda_3\le 1/n$.
In Section 8, the main theorem~2 is proved, which gives an upper bound for $\Lambda_n^2$.
Everywhere we assume that $\sum\limits_{k=s}^pa_k=0$ and $\prod\limits_{k=s}^pa_k=1$ if $p<s$.

\section{Lagrange's Fundamental Interpolating Polynomials}

It is known that
$$
l_i(\lambda)=\frac{\omega_i(\lambda)}{\omega_i(a_i)}~,\quad i\in I,
$$
where
$$
\omega_{i}(\lambda)=
\prod\limits_{s_1=0}^{i_1-1}\left(\lambda_1-\frac{s_1}{n}\right)
\prod\limits_{s_2=0}^{i_2-1}\left(\lambda_2-\frac{s_2}{n}\right)
\ldots
\prod\limits_{s_{d+1}=0}^{i_{d+1}-1}\left(\lambda_{d+1}-\frac{s_{d+1}}{n}
\right)={}
$$
$$
{}=\left(\frac{1}{n}\right)^n
\frac{\Gamma(n\lambda_1+1)}{\Gamma(n\lambda_1-i_1+1)}
~\ldots~
\frac{\Gamma(n\lambda_{d+1}+1)}{\Gamma(n\lambda_{d+1}-i_{d+1}+1)}~.
$$
This explicit formula was given by Nicolaides~\cite{Niko1972}.
Indeed,
$\omega_i(a_j)=0$ for $i\neq j$ since if $i,j\in I$,
$i\neq j$, then there is
$q\in \{1,\dots,d+1\}$ such than $j_q<i_q$ and hence
$
{\displaystyle
\prod\limits_{s_q=0}^{i_q-1}
\left(\frac{j_q}{n}-\frac{s_q}{n}\right)=0.
}
$
Since
$
\omega_i(a_i)\!=\!
\left( 1/n \right)^n\Gamma(i_1\!+\!1)\ldots\Gamma(i_{d+1}\!+\!1)=
\left( 1/n \right)^ni_1!i_2!\ldots i_{d+1}!,
$
then
\begin{equation}
\label{fund_pol}
l_i(\lambda)=
\frac{\Gamma(n\lambda_1+1)\ldots\Gamma(n\lambda_{d+1}+1)}
     {\Gamma(n\lambda_1-i_1+1)\ldots\Gamma(n\lambda_{d+1}-i_{d+1}+1)}
     \cdot
     \frac{1}{i_1!i_2!\ldots i_{d+1}!}~.
\end{equation}
When obtaining estimates, we will exclude from
consideration those arguments of $\Gamma$-functions that
are poles, since these cases correspond to the situation
when the corresponding summands $|l_i(\lambda)|$
from (\ref{funcLeb}) are equal to zero.

%When obtaining estimates, we will exclude from consideration
%those $\lambda$ and $i$ for which the argument of at least
%one of the $\Gamma$--functions of the denominator is its pole,
%since this case corresponds to the situation $l_i(\lambda) =0,$ i.e.,
%the corresponding summand in (\ref{funcLeb}) is equal to zero.

\section{Properties of the Lebesgue Function}

Further, let $d=2$, $\Delta$ is a triangle, $n\ge 4$,
${\mathcal{L}}_n={\mathcal{L}}_n(\lambda)={\mathcal{L}}_n^2(\lambda)$.

\begin{lemen}
\label{symm}
The Lebesgue function ${\mathcal{L}}_n(\lambda_1,\lambda_2,\lambda_3)$
is a symmetric function of its arguments.
\end{lemen}

\textsl{Proof}.
Let's prove, for example, that
${\mathcal{L}}_n(\lambda_1,\lambda_2,\lambda_3)={\mathcal{L}}_n(\lambda_2,\lambda_1,\lambda_3)$.
This follows from~(\ref{funcLeb}) and~(\ref{fund_pol}),
because for any $i=(i_1,i_2,i_3)\in I$
there is $j=(j_1,j_2,j_3) =(i_2,i_1,i_3)\in I$
for which
$
l_i(\lambda_1,\lambda_2,\lambda_3)
=
l_j(\lambda_2,\lambda_1,\lambda_3)$.
The other cases are proved in the same way.
\hfill$\Box$

\begin{lemen}
\label{restrict}
%Справедливо равенство
$
%\begin{equation}
%\label{restrict_form}
\quad
\Lambda_n
=
\underset
{
\lambda_1\ge\lambda_2\ge 0
\atop
{
\lambda_1\ge\lambda_3\ge 0
\atop
0\le\lambda_1+\lambda_2+\lambda_3\le1
}
}
\max
{\mathcal{L}}_n(\lambda_1,\lambda_2,\lambda_3).
%\end{equation}
$
\end{lemen}

\textsl{Proof}.
%\begin{proof}
This equality follows from the fact that
the function ${\mathcal{L}}_n(\lambda_1,\lambda_2,\lambda_3)$
is symmetric.
\hfill$\Box$

\vspace{5mm}

Note that if a point $u$ with barycentric
coordinates $(\lambda_1,\lambda_2,\lambda_3)$ belongs to $\Delta$
then there are numbers $r_1,r_2,r_3\in{\Bbb Z}_+$,
$r_1+r_2+r_3=n-1$,
$\alpha_1\in[-1,1],$
$\alpha_2,\alpha_3\in[0,1]$,
$\alpha_1+\alpha_2+\alpha_3=1$,
such that $n\lambda_s=r_s+\alpha_s$, $s=1,2,3$
(see, for example, Fig.~1).

\begin{center}
\includegraphics[height=6cm]{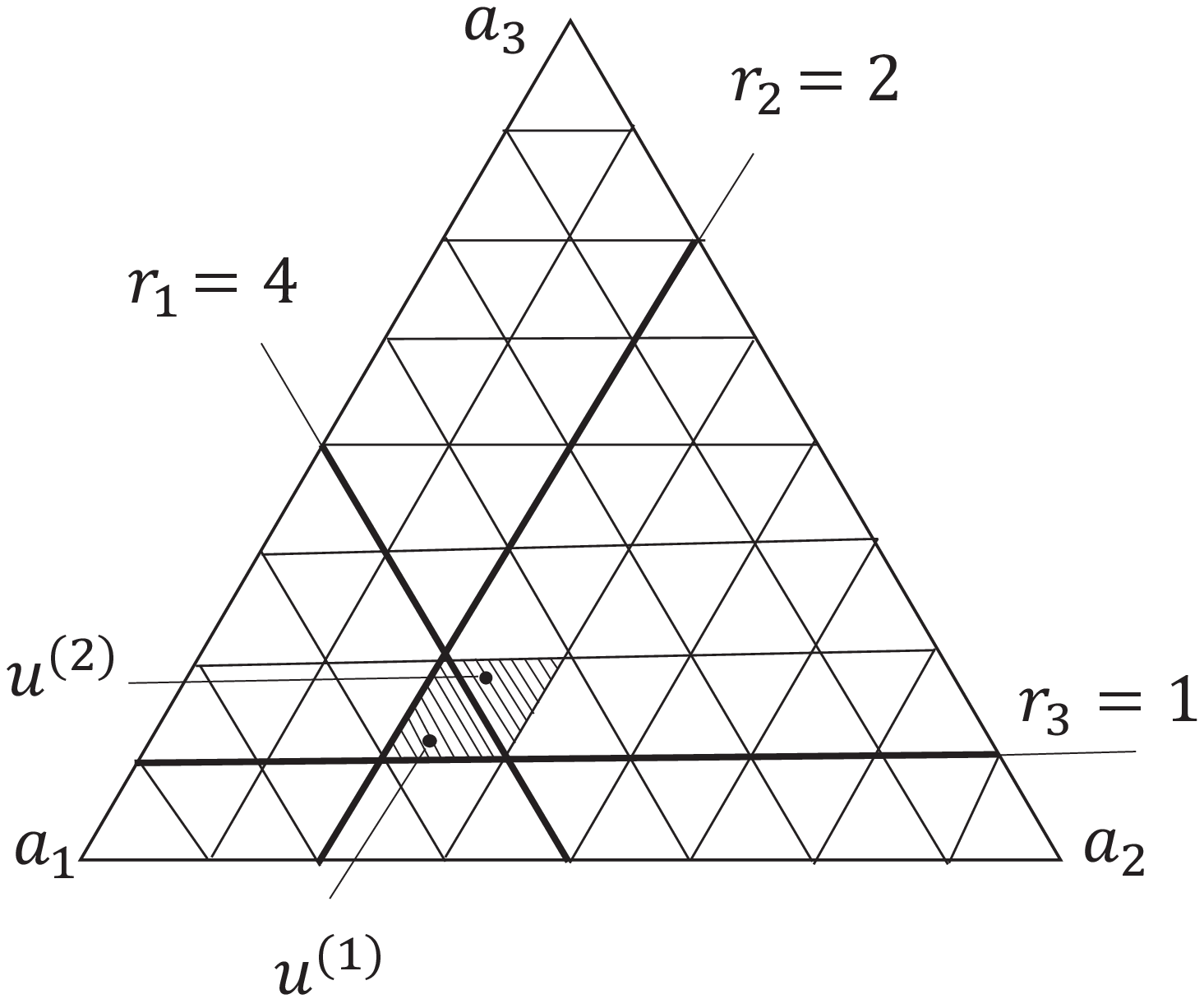}
\end{center}

{\bf Fig. 1.} {\small Equally spaced nodes in the case $n=8$. Here
$r_1=4$, $r_2=2$, $r_3=1$,
$0<\alpha_1, \alpha_2, \alpha_3<1$ for $u^{(1)}$
or
$r_1=4$, $r_2=2$, $r_3=1$,
$-1<\alpha_1<0$, $0<\alpha_2, \alpha_3<1$ for $u^{(2)}$.}

\vspace{5mm}

So, let $n\lambda_s=r_s+\alpha_s$, $r_1+r_2+r_3=n-1$, $-1\le \alpha_1\le 1$,
$0\le \alpha_2,\alpha_3\le 1$. Denote
$
\mathfrak{L}_n(r_1,r_2,r_3,\alpha_1,\alpha_2,\alpha_3)
=
{\mathcal{L}}_n(\lambda_1,\lambda_2,\lambda_3),
$
and
$$
a_{i_s}(\lambda_s)=
\left|
\frac{\Gamma(n\lambda_s+1)}
     {i_s!\Gamma(n\lambda_s-i_s+1)}
\right|=
\frac{\Gamma(r_s+\alpha_s+1)}{i_s!\left|\Gamma(r_s+\alpha_s-i_s+1)\right|}~,
\quad s=1,2,3.
$$
Also, let us define the sums
$S_k=S_k(r_1,r_2,r_3,\alpha_1,\alpha_2,\alpha_3)$,
$1\le k\le 6$, as follows:
$$
S_1=
\sum\limits_{i_2=0}^{r_2}
\sum\limits_{i_3=0}^{r_3}
|l_i(\lambda)|,
\quad
S_3=
\sum\limits_{i_2=0}^{r_2}
\sum\limits_{i_3=0}^{r_1}
|l_i(\lambda)|,
$$
$$
S_4=
\sum\limits_{i_2=0}^{r_2-1}
\sum\limits_{i_3=r_3+1}^{n-r_1-1-i_2}
|l_i(\lambda)|,
\quad
S_5=
\sum\limits_{i_2=r_2+1}^{n-r_1-1}
\sum\limits_{i_3=0}^{n-i_2}
|l_i(\lambda)|,
\quad
S_6=
\sum\limits_{i_2=r_2+1}^{n-r_3-1}
\sum\limits_{i_1=0}^{n-r_3-1-i_2}
|l_i(\lambda)|.
$$

$$
S_2=S_{2,1}+S_{2,2}+S_{2,3},
$$
where
$$
S_{2,1}=
\sum\limits_{i_2=n-r_1}^{n-r_3}
\sum\limits_{i_1=n-r_3-i_2}^{n-i_2}
|l_i(\lambda)|,
~
S_{2,2}=
\sum\limits_{i_2=n-r_3+1}^{n}
\sum\limits_{i_1=0}^{n-i_2}
|l_i(\lambda)|,
~
S_{2,3}=
\sum\limits_{i_2=r_2+1}^{n-r_1-1}
\sum\limits_{i_1=n-r_3-i_2}^{r_1}
|l_i(\lambda)|.
$$
Then
\begin{equation}
\label{abs_l}
|l_i(\lambda)|
=
\prod\limits_{s=1}^3
a_{i_s}(\lambda_s)
=
\prod\limits_{s=1}^3
\frac{\Gamma(r_s+\alpha_s+1)}{i_s!\left|\Gamma(r_s+\alpha_s-i_s+1)\right|}~,
\end{equation}

\begin{equation}
\label{Lambda_S}
\mathfrak{L}_n(r_1,r_2,r_3,\alpha_1,\alpha_2,\alpha_3)
=
{\mathcal{L}}_n(\lambda_1,\lambda_2,\lambda_3)
=
\sum\limits_{k=1}^{6}S_k.
%\overset{def}{=}
%{\mathcal S}(r_1,r_2,r_3,\alpha_1,\alpha_2,\alpha_3),
\end{equation}

%==========================================================================

Taking into account~(\ref{const_Leb}) and Lemma~\ref{restrict},
we can state that in order to estimate $\Lambda_n$ it is sufficient to estimate the maximum of the function $\mathfrak{L}_n(r_1,r_2,r_3,\alpha_1,\alpha_2,\alpha_3)
= \mathcal{L}_n(\lambda_1,\lambda_2,\lambda)$
for
$\lambda_2\le\lambda_1$, $\lambda_3\le \lambda_1$,
so we further assume that
\begin{equation}
\label{constr_arg}
r_1+\alpha_1\ge r_2+\alpha_2,
\quad
r_1+\alpha_1\ge r_3+\alpha_3
\end{equation}
and we will estimate
$\mathfrak{L}_n=\mathcal{L}_n$
under constraints~(\ref{constr_arg}).
Denote
\begin{equation}
\label{delta_k}
\delta_k={S}_k(r_1,r_2,r_3,\alpha_1,\alpha_2,\alpha_3)-{S}_k(r_1-1,r_2,r_3+1,\alpha_1,\alpha_2,\alpha_3),
\quad
1\le k\le 6.
\end{equation}

The purpose of Section 5 is to estimate the values of
$\delta_k$
under the constraints
$r_1-1+\alpha_1\ge r_2+\alpha_2,$
$r_1-1+\alpha_1\ge r_3+1+\alpha_3$,
$-1<\alpha_1<1.$
Since $r_1,r_2,r_3\in\mathbb{Z}_+$ then
$$
r_1\ge r_3+2\ge 2.
$$

%==========================================

\section{Some Auxiliary Results}

We will use the notation
$
\displaystyle{
\left(
\begin{array}{c}
a\\
b
\end{array}
\right)
=
\frac{\Gamma(a+1)}{\Gamma(b+1)\Gamma(a-b+1)}~.
}
$

\begin{lemen}
\label{comb}
%Имеет место равенство
$
\quad
\displaystyle{
%\begin{equation}
%\label{gamma_sum}
\frac{1}
{b
\left(
\begin{array}{c}
a-1\\
b
\end{array}
\right)}
-
\frac{1}
{(b+1)
\left(
\begin{array}{c}
a\\
b+1
\end{array}
\right)}
=
\frac{1}
{b
\left(
\begin{array}{c}
a\\
b
\end{array}
\right)}~.
%\end{equation}
}
$
\end{lemen}

\textsl{Proof}.
The proof
follows from the fact that
$$
\frac{1}
{(b+1)
\left(
\begin{array}{c}
a\\
b+1
\end{array}
\right)}
+
\frac{1}
{b
\left(
\begin{array}{c}
a\\
b
\end{array}
\right)}
=
\frac{\Gamma(b)\Gamma(a-b)}
{\Gamma(a+1)}
(b+a-b)
=
\frac{1}
{b
\left(
\begin{array}{c}
a-1\\
b
\end{array}
\right)}~.
$$
\hfill$\Box$

\begin{lemen}
\label{vspom}
Let
$
\displaystyle{
{\mathcal D}^*=
{\mathcal D}^*(p,m,x,y)=
\sum\limits_{\kappa=0}^{p+1}
\frac{
\left(
\begin{array}{c}
x+1\\
\kappa
\end{array}
\right)
}
{y
\left(
\begin{array}{c}
m-\kappa\\
y
\end{array}
\right)
}
-
\sum\limits_{\kappa=0}^{p}
\frac{
\left(
\begin{array}{c}
x\\
\kappa
\end{array}
\right)
}
{(y+1)
\left(
\begin{array}{c}
m-\kappa\\
y+1
\end{array}
\right)
}~.
%\right)
}
$
%где $r,q,m,p\in{\Bbb Z}_+$, %$p\ge r$, $-1\le \alpha,\beta\le 1$.
Then
\begin{equation}
\label{lemen_vspom}
{\mathcal D}^*=
\sum\limits_{\kappa=0}^{p}
\frac{2
\left(
\begin{array}{c}
x\\
\kappa
\end{array}
\right)
}
{y
\left(
\begin{array}{c}
m-\kappa\\
y
\end{array}
\right)
}
+
\frac{
\left(
\begin{array}{c}
x\\
p+1
\end{array}
\right)
}
{y
\left(
\begin{array}{c}
m-p-1\\
y
\end{array}
\right)
}~.
\end{equation}
\end{lemen}

\textsl{Proof}.
Let's make the following transformations:
$$
{\mathcal D}^*=
%\frac{
%\left(
%\begin{array}{c}
%r+1+\alpha\\
%{p+1}
%\end{array}
%\right)
%}
%{(q+\beta)
%\left(
%\begin{array}{c}
%m-p-1\\
%q+\beta
%\end{array}
%\right)
%}
%+
\sum\limits_{\kappa=0}^{p+1}
\frac{
\left(
\begin{array}{c}
x\\
\kappa
\end{array}
\right)
}
{y
\left(
\begin{array}{c}
m-\kappa\\
y
\end{array}
\right)
}
+
\sum\limits_{\kappa=1}^{p+1}
\frac{
\left(
\begin{array}{c}
x\\
\kappa-1
\end{array}
\right)
}
{y
\left(
\begin{array}{c}
m-\kappa\\
y
\end{array}
\right)
}
-
\sum\limits_{\kappa=0}^{p}
\frac{
\left(
\begin{array}{c}
x\\
\kappa
\end{array}
\right)
}
{(y+1)
\left(
\begin{array}{c}
m-\kappa\\
y+1
\end{array}
\right)
}~.
$$
Let's change the variable $\kappa=\nu+1$
in the second sum and use the lemma~\ref{comb}
for $a=m-\nu-1$, $b=y$.
Then
$$
{\mathcal D}^*
=
\sum\limits_{\kappa=0}^{p+1}
\frac{
\left(
\begin{array}{c}
x\\
\kappa
\end{array}
\right)
}
{y
\left(
\begin{array}{c}
m-\kappa\\
y
\end{array}
\right)
}
+
\left(
\sum\limits_{\nu=0}^{p}
\frac{
\left(
\begin{array}{c}
x\\
\nu
\end{array}
\right)
}
{y
\left(
\begin{array}{c}
m-\nu-1\\
y
\end{array}
\right)
}
{}
-
\sum\limits_{\kappa=0}^{p}
\frac{
\left(
\begin{array}{c}
x\\
\kappa
\end{array}
\right)
}
{(y+1)
\left(
\begin{array}{c}
m-\kappa\\
y+1
\end{array}
\right)
}
\right)
$$
$$
{}
=
\sum\limits_{\kappa=0}^{p+1}
\frac{
\left(
\begin{array}{c}
x\\
\kappa
\end{array}
\right)
}
{y
\left(
\begin{array}{c}
m-\kappa\\
y
\end{array}
\right)
}
+
\sum\limits_{\kappa=0}^{p}
\frac{
\left(
\begin{array}{c}
x\\
\kappa
\end{array}
\right)
}
{y
\left(
\begin{array}{c}
m-\kappa\\
y
\end{array}
\right)
}~,
$$
and~(\ref{lemen_vspom}) follows from this equality.
\hfill$\Box$

%===============================================================

%======================================================
\begin{lemen}
\label{vspom1}
Let
$$
{\mathcal D}^{**}=
{\mathcal D}^{**}(m,q,r,x,y)=
\sum\limits_{\kappa=q}^{m-1-r}
\frac{1}
{
(x+1)
\left(
\begin{array}{c}
m-\kappa\\
x+1
\end{array}
\right)
y
\left(
\begin{array}{c}
\kappa\\
y
\end{array}
\right)
}
-
\sum\limits_{\kappa=q+1}^{m-r}
\frac{1}
{
x
\left(
\begin{array}{c}
m-\kappa\\
x
\end{array}
\right)
(y+1)
\left(
\begin{array}{c}
\kappa\\
y+1
\end{array}
\right)
}~.
$$
Then
$$
{\mathcal D}^{**}
=
\frac{1}
{
x
\left(
\begin{array}{c}
r\\
x
\end{array}
\right)
y
\left(
\begin{array}{c}
m-r\\
y
\end{array}
\right)
}
-
\frac{1}
{
x
\left(
\begin{array}{c}
m-q\\
x
\end{array}
\right)
y
\left(
\begin{array}{c}
q\\
y
\end{array}
\right)
}~.
%\end{equation}
$$
\end{lemen}

\textsl{Proof}.
Using Lemma~\ref{comb} with $a=m-\kappa$, $b=r+\alpha$ we have
$$
{\mathcal D}^{**}
=
\sum\limits_{\kappa=q}^{m-1-r}
\frac{1}
{
x
\left(
\begin{array}{c}
m-\kappa-1\\
x
\end{array}
\right)
y
\left(
\begin{array}{c}
\kappa\\
y
\end{array}
\right)
}
-
\sum\limits_{\kappa=q}^{m-1-r}
\frac{1}
{
x
\left(
\begin{array}{c}
m-\kappa\\
x
\end{array}
\right)
y
\left(
\begin{array}{c}
\kappa\\
y
\end{array}
\right)
}
$$
$$
{}
-
\sum\limits_{\kappa=q+1}^{m-r}
\frac{1}
{
x
\left(
\begin{array}{c}
m-\kappa\\
x
\end{array}
\right)
(y+1)
\left(
\begin{array}{c}
\kappa\\
y+1
\end{array}
\right)
}~.
$$
In the first sum, we make the change of variable $\kappa=\nu-1$.
Add terms with the same $\kappa$ from the second and third sums,
leaving the term with $\kappa=q$ in the second sum and the term
with $\kappa=m-r$ in the third sum unchanged.
Then
$$
{\mathcal D}^{**}
=
\sum\limits_{\nu=q+1}^{m-r}
\frac{1}
{
x
\left(
\begin{array}{c}
m-\nu\\
x
\end{array}
\right)
y
\left(
\begin{array}{c}
\nu-1\\
y
\end{array}
\right)
}
-
\sum\limits_{\kappa=q+1}^{m-1-r}
\frac{1}
{
x
\left(
\begin{array}{c}
m-\kappa\\
x
\end{array}
\right)
y
\left(
\begin{array}{c}
\kappa-1\\
y
\end{array}
\right)
}
$$
$$
{}
-
\frac{1}
{
x
\left(
\begin{array}{c}
m-q\\
x
\end{array}
\right)
y
\left(
\begin{array}{c}
q\\
y
\end{array}
\right)
}
-
\frac{1}
{
x
\left(
\begin{array}{c}
r\\
x
\end{array}
\right)
(y+1)
\left(
\begin{array}{c}
m-r\\
y+1
\end{array}
\right)
}
$$
$$
{}=
\left(
\frac{1}
{
x
\left(
\begin{array}{c}
r\\
x
\end{array}
\right)
y
\left(
\begin{array}{c}
m-r-1\\
y
\end{array}
\right)
}
-
\frac{1}
{
x
\left(
\begin{array}{c}
r\\
x
\end{array}
\right)
(y+1)
\left(
\begin{array}{c}
m-r\\
y+1
\end{array}
\right)
}
\right)
-
\frac{1}
{
x
\left(
\begin{array}{c}
m-q\\
x
\end{array}
\right)
y
\left(
\begin{array}{c}
q\\
y
\end{array}
\right)
}
$$
$$
{}
=
\frac{1}
{
x
\left(
\begin{array}{c}
r\\
x
\end{array}
\right)
x
\left(
\begin{array}{c}
m-r\\
x
\end{array}
\right)
}
-
\frac{1}
{
x
\left(
\begin{array}{c}
m-q\\
x
\end{array}
\right)
y
\left(
\begin{array}{c}
q+1\\
y
\end{array}
\right)
}~.
$$
\hfill$\Box$

%======================================================
\begin{lemen}
\label{vspom2}
Let
$
\displaystyle{
{\mathcal D}^{***}=
{\mathcal D}^{***}(p,q,m,x,y)=
\sum\limits_{\kappa=q}^{p}
\left(
\begin{array}{c}
x\\
\kappa
\end{array}
\right)
\left(
\begin{array}{c}
y\\
m-\kappa
\end{array}
\right)
-
\sum\limits_{\kappa=q-1}^{p-1}
{
\left(
\begin{array}{c}
x-1\\
\kappa
\end{array}
\right)
\left(
\begin{array}{c}
y+1\\
m-\kappa
\end{array}
\right)
}.
}
$
Then
$$
{\mathcal D}^{***}=
\left(
\begin{array}{c}
x-1\\
p
\end{array}
\right)
\left(
\begin{array}{c}
y\\
m-p
\end{array}
\right)
-
{
\left(
\begin{array}{c}
x-1\\
q-1
\end{array}
\right)
\left(
\begin{array}{c}
y\\
m-q+1
\end{array}
\right).
}~
$$
\end{lemen}

\textsl{Proof}. Notice that
$$
D^{***}=
\sum\limits_{\kappa=q}^{p-1}
\left(
\begin{array}{c}
x-1\\
\kappa
\end{array}
\right)
\left(
\begin{array}{c}
y\\
m-\kappa
\end{array}
\right)
+
\sum\limits_{\kappa=q}^{p-1}
\left(
\begin{array}{c}
x-1\\
\kappa-1
\end{array}
\right)
\left(
\begin{array}{c}
y\\
m-\kappa
\end{array}
\right)
+
\left(
\begin{array}{c}
x\\
p
\end{array}
\right)
\left(
\begin{array}{c}
y\\
m-p
\end{array}
\right)
$$
$$
{}
-
\sum\limits_{\kappa=q}^{p-1}
\left(
\begin{array}{c}
x-1\\
\kappa
\end{array}
\right)
\left(
\begin{array}{c}
y\\
m-\kappa
\end{array}
\right)
-
\sum\limits_{\kappa=q}^{p-1}
\left(
\begin{array}{c}
x-1\\
\kappa
\end{array}
\right)
\left(
\begin{array}{c}
y\\
m-\kappa-1
\end{array}
\right)
-
\left(
\begin{array}{c}
x-1\\
q-1
\end{array}
\right)
\left(
\begin{array}{c}
y+1\\
m-q+1
\end{array}
\right).
$$
The first sum in the first line and the first
sum in the second line give $0$, in the second
sum of the first line we make the change
of variable $\kappa=s+1$. Then$$
{\mathcal D}^{***}=
\sum\limits_{s=q-1}^{p-2}
\left(
\begin{array}{c}
x-1\\
s
\end{array}
\right)
\left(
\begin{array}{c}
y\\
m-s-1
\end{array}
\right)
-
\sum\limits_{\kappa=q}^{p-1}
\left(
\begin{array}{c}
x-1\\
\kappa
\end{array}
\right)
\left(
\begin{array}{c}
y\\
m-\kappa-1
\end{array}
\right)
+
\left(
\begin{array}{c}
x-1\\
p
\end{array}
\right)
\left(
\begin{array}{c}
y\\
m-p
\end{array}
\right)
$$
$$
{}
+
\left(
\begin{array}{c}
x-1\\
p-1
\end{array}
\right)
\left(
\begin{array}{c}
y\\
m-p
\end{array}
\right)
-
\left(
\begin{array}{c}
x-1\\
q-1
\end{array}
\right)
\left(
\begin{array}{c}
y\\
m-q+1
\end{array}
\right)
-
\left(
\begin{array}{c}
x-1\\
q-1
\end{array}
\right)
\left(
\begin{array}{c}
y\\
m-q
\end{array}
\right)
$$

$$
{}=
\sum\limits_{s=q-1}^{p-1}
\left(
\begin{array}{c}
x-1\\
s
\end{array}
\right)
\left(
\begin{array}{c}
y\\
m-s-1
\end{array}
\right)
-
\sum\limits_{\kappa=q-1}^{p-1}
\left(
\begin{array}{c}
x-1\\
\kappa
\end{array}
\right)
\left(
\begin{array}{c}
y\\
m-\kappa-1
\end{array}
\right)
\left(
\begin{array}{c}
x-1\\
p
\end{array}
\right)
\left(
\begin{array}{c}
y\\
m-p
\end{array}
\right)
$$
$$
{}
-
\left(
\begin{array}{c}
x-1\\
q-1
\end{array}
\right)
\left(
\begin{array}{c}
y\\
m-q+1
\end{array}
\right)
=
\left(
\begin{array}{c}
x-1\\
p
\end{array}
\right)
\left(
\begin{array}{c}
y\\
m-p
\end{array}
\right)
-
\left(
\begin{array}{c}
x-1\\
q-1
\end{array}
\right)
\left(
\begin{array}{c}
y\\
m-q+1
\end{array}
\right).
$$
\hfill$\Box$
%====================================================

\begin{lemen}
\label{monot}
Let
$
\displaystyle{
g(\alpha)
=
\frac{\Gamma(a+\alpha)}{\Gamma(b+\alpha)}.
}
$
If $a\ge b\ge \alpha$, then the function
$g(\alpha)$ is nondecreasing.
If $b\ge a\ge \alpha$, then the function
$g(\alpha)$
is nonincreasing.
\end{lemen}

\textsl{Proof}.
The assertion of the lemma follows from the fact that
$$
g'(\alpha)
=
\frac{\Gamma(a+\alpha)\psi(a+\alpha)\Gamma(b+\alpha)-\Gamma(a+\alpha)\Gamma(b+\alpha)\psi(b+\alpha)}{(\Gamma(b+\alpha))^2}
$$
$$
{}
=
\frac{\Gamma(a+\alpha)\Gamma(b+\alpha)}{(\Gamma(b+\alpha))^2}
\left(\psi(a+\alpha)-\psi(b+\alpha)\right),
$$
where $\psi(t)$ is the digamma function.
Since $\psi(t)$ increases on $(0,+\infty)$
then $g'(\alpha)\ge 0$ for $a\ge b$ and $g'(\alpha)\le 0$ for $ a\le b$.
\hfill$\Box$

\section{Estimations of $\delta_s$, $1\le \delta_s\le 6$}
Recall that we exclude from consideration those
$\lambda$ and $i$
for which the argument of at least one of the
$\Gamma$--functions in the denominators
is its pole, since these cases correspond
to the situation $l_i(\lambda )=0,$ i.e.,
the corresponding term in (\ref{funcLeb}) is equal to zero.
Also recall that we are considering
the case
$r_1-1+\alpha_1\ge r_2+\alpha_2$,
$r_1-1+\alpha_1\ge r_3+1+\alpha_3$.
In particular, $r_1\ge r_3+2\ge 2$.
Everywhere below we assume that
$i=(i_1,i_2,i_3)\in I$, i.e. $i_1+i_2+i_3=n.$

\begin{lemen}
\label{delta_1}
$
\displaystyle{
\delta_1
\ge
-2^{r_2+r_3+2}+2^{r_3+1}.
}
$
\end{lemen}

\textsl{Proof}.
Due to Lemma~\ref{vspom}
with
$p=r_3$, $m=n-i_2$, $x=r_3+\alpha_3$, $y=r_1+\alpha_1$ we have
$$
\delta_1=
\frac{|\sin \pi\alpha_1|}{\pi}
\sum\limits_{i_2=0}^{r_2} a_{i_2}(\lambda_2)
\sum\limits_{i_3=0}^{r_3}
\frac{\Gamma(r_1+\alpha_1+1)\Gamma(i_1-r_1-\alpha_1)}
{i_1!}
~
\frac{\Gamma(r_3+\alpha_3+1)}
{i_3!~\Gamma(r_3+\alpha_3+1-i_3)}
$$
$$
{}-
\frac{|\sin \pi\alpha_1|}{\pi}
\sum\limits_{i_2=0}^{r_2} a_{i_2}(\lambda_2)
\sum\limits_{i_3=0}^{r_3+1}
\frac{\Gamma(r_1+\alpha_1)\Gamma(i_1-r_1-\alpha_1+1)}
{i_1!}
~
\frac{\Gamma(r_3+\alpha_3+2)}
{i_3!~\Gamma(r_3+\alpha_3+2-i_3)}
$$
$$
{}=-
\frac{|\sin \pi\alpha_1|}{\pi}
\sum\limits_{i_2=0}^{r_2} a_{i_2}(\lambda_2)
\left(
\sum\limits_{i_3=0}^{r_3+1}
\frac{
\left(
\begin{array}{c}
r_3+1+\alpha_3\\
i_3
\end{array}
\right)
}
{(r_1+\alpha_1)
\left(
\begin{array}{c}
i_1\\
r_1+\alpha_1
\end{array}
\right)
}
-
\sum\limits_{i_3=0}^{r_3}
\frac{
\left(
\begin{array}{c}
r_3+\alpha_3\\
i_3
\end{array}
\right)
}
{(r_1+\alpha_1+1)
\left(
\begin{array}{c}
i_1\\
r_1+\alpha_1+1
\end{array}
\right)
}
\right)
$$
$$
{}
=
-
\frac{|\sin \pi\alpha_1|}{\pi}
\sum\limits_{i_2=0}^{r_2} a_{i_2}(\lambda_2)
{\mathcal D}^*(r_3,n-i_2, r_3+\alpha_3,r_1+\alpha_1)
$$
$$
{}
=
-
\frac{|\sin \pi\alpha_1|}{\pi}
\sum\limits_{i_2=0}^{r_2} a_{i_2}(\lambda_2)
\left(
\sum\limits_{i_3=0}^{r_3}
\frac{2
\left(
\begin{array}{c}
r_3+\alpha_3\\
i_3
\end{array}
\right)
}
{(r_1+\alpha_1)
\left(
\begin{array}{c}
n-i_2-i_3\\
r_1+\alpha_1
\end{array}
\right)
}
+
\frac{
\left(
\begin{array}{c}
r_3+\alpha_3\\
r_3+1
\end{array}
\right)
}
{(r_1+\alpha_1)
\left(
\begin{array}{c}
n-i_2-r_3-1\\
r_1+\alpha_1
\end{array}
\right)
}
\right).
$$
Since $0\le i_2\le r_2$, $0\le i_3\le r_3$, $i_1=n-i_2-i_3\ge n-r_2-r_3=r_1+1$,
then all the terms in the obtained sums are non-negative
and all the gamma functions in fractions
$
\displaystyle{
\frac{\Gamma(a+1)}{\Gamma(b+1)\Gamma(a-b+1)}
=
\left(
\begin{array}{c}
a\\b
\end{array}
\right)
}
$
are non-negative.
Then
$$
\delta_1
\ge -
\frac{|\sin \pi\alpha_1|}{\pi}
\sum\limits_{i_2=0}^{r_2} a_{i_2}(\lambda_2)
\left(
\sum\limits_{i_3=0}^{r_3}
\frac{2
\left(
\begin{array}{c}
r_3+1\\
i_3
\end{array}
\right)
}
{(r_1+\alpha_1)
\left(
\begin{array}{c}
r_1+1\\
r_1+\alpha_1
\end{array}
\right)
}
+
\frac{1}
{(r_1+\alpha_1)
\left(
\begin{array}{c}
r_1\\
r_1+\alpha_1
\end{array}
\right)
}
\right)
$$
$$
{}=-
\frac{|\sin \pi\alpha_1|}{\pi}
\sum\limits_{i_2=0}^{r_2} a_{i_2}(\lambda_2)
\left(
\sum\limits_{i_3=0}^{r_3}
\frac{2\Gamma(r_1+\alpha_1)\Gamma(2-\alpha_1)}{(r_1+1)!}
\left(
\begin{array}{c}
r_3+1\\
i_3
\end{array}
\right)
+
\frac{\Gamma(r_1+\alpha_1)\Gamma(1-\alpha_1)}
{r_1!}
\right)
$$
$$
{}\ge -
\sum\limits_{i_2=0}^{r_2} a_{i_2}(\lambda_2)
\sum\limits_{i_3=0}^{r_3}
\frac{2(1-\alpha_1)}{(r_1+1)|\Gamma(\alpha_1)|}
\left(
\begin{array}{c}
r_3+1\\
i_3
\end{array}
\right)
-
\sum\limits_{i_2=0}^{r_2}
\frac{a_{i_2}(\lambda_2)}
{|\Gamma(\alpha_1)|}
\ge -
\sum\limits_{i_2=0}^{r_2} a_{i_2}(\lambda_2)
\sum\limits_{i_3=0}^{r_3+1}
\left(
\begin{array}{c}
r_3+1\\
i_3
\end{array}
\right)
$$
$$
{}=-
\sum\limits_{i_2=0}^{r_2}
\frac{\Gamma(r_2+\alpha_2+1)}{i_2!~\Gamma(r_2+\alpha_2+1-i_2)}
2^{r_3+1} %%+\pi-1
\ge -
2^{r_3+1}
\sum\limits_{i_2=0}^{r_2}
\frac{(r_2+1)!}{i_2!~(r_2+1-i_2)!}
\ge
%-2^{r_3+1}(2^{r_2+1}-1)=
-2^{r_2+r_3+2}+2^{r_3+1}
%=2^{n+1-r_1},
$$
(we used Lemma~\ref{monot} in the last but one inequality).
\hfill$\Box$

\begin{lemen}
\label{delta_3}
%Имеет место неравенство
$
\delta_3\ge -\frac{2^{r_2-1}}{r_1}.
$
%\begin{equation}
%\label{lemen_delta_3}
%\delta_3=
%5\frac{\sin \pi\alpha_3}{\pi}
%\sum\limits_{i_2=0}^{r_2} a_{i_2}
%\left(
%\sum\limits_{i_1=0}^{r_1-1}
%\frac{2
%\left(
%\begin{array}{c}
%r_1-1+\alpha_1\\
%i_1
%\end{array}
%\right)
%}
%{(r_3+\alpha_3+1)
%\left(
%\begin{array}{c}
%i_3\\
%r_3+\alpha_3+1
%\end{array}
%\right)
%}
%+
%\frac{
%\left(
%\begin{array}{c}
%r_1-1+\alpha_1\\
%r_1
%\end{array}
%\right)
%}
%{(r_3+\alpha_3+1)
%\left(
%\begin{array}{c}
%n-i_2-r_1\\
%r_3+\alpha_3+1
%\end{array}
%\right)
%}
%\right).
%\end{equation}
\end{lemen}

\textsl{Proof}.
Due to Lemma~\ref{vspom}
with
$p-r_1-1$, $m-n-i_2$, $x=r_1+\alpha_1-1$, $y=r_3+\alpha_3+1$
we obtain
$$
\delta_3=
\frac{\sin \pi\alpha_3}{\pi}
\sum\limits_{i_2=0}^{r_2} a_{i_2}(\lambda_2)
\sum\limits_{i_1=0}^{r_1}
\frac{\Gamma(r_3+\alpha_3+1)\Gamma(i_3-r_3-\alpha_3)}{i_3!}
%\left|
%\frac{\Gamma(r_3+\alpha_3+1)}
%{i_3!~\Gamma(r_3+\alpha_3+1-i_3)}
%\right|
~
\frac{\Gamma(r_1+\alpha_1+1)}
{i_1!~\Gamma(r_1+\alpha_1+1-i_1)}
$$
$$
{}-
\frac{\sin \pi\alpha_3}{\pi}
\sum\limits_{i_2=0}^{r_2} a_{i_2}(\lambda_2)
\sum\limits_{i_1=0}^{r_1-1}
\frac{\Gamma(r_3+\alpha_3+2)\Gamma(i_3-r_3-1-\alpha_3)}
{i_3!}
~
\frac{\Gamma(r_1+\alpha_1)}
{i_1!~\Gamma(r_1+\alpha_1-i_1)}
$$

$$
{}=
\frac{\sin \pi\alpha_3}{\pi}
\sum\limits_{i_2=0}^{r_2} a_{i_2}(\lambda_2)
{\mathcal D}^*(r_1-1,n-i_2, r_1+\alpha_1-1,r_3+\alpha_3+1)
$$
$$
{}
=
\frac{\sin \pi\alpha_3}{\pi}
\sum\limits_{i_2=0}^{r_2} a_{i_2}(\lambda_2)
\left(
\sum\limits_{i_1=0}^{r_1-1}
\frac{2
\left(
\begin{array}{c}
r_1-1+\alpha_1\\
i_1
\end{array}
\right)
}
{(r_3+\alpha_3+1)
\left(
\begin{array}{c}
n-i_2-i_1\\
r_3+\alpha_3+1
\end{array}
\right)
}
+
\frac{
\left(
\begin{array}{c}
r_1-1+\alpha_1\\
r_1
\end{array}
\right)
}
{(r_3+\alpha_3+1)
\left(
\begin{array}{c}
n-i_2-r_1\\
r_3+\alpha_3+1
\end{array}
\right)
}
\right).
%\ge 0
$$

Let $\alpha_1\ge 0$.
Since $0\le i_1\le r_1-1$, $0\le i_2\le r_2$,
$i_3=n-i_1-i_2\ge n-r_1+1-r_2=r_3$,
then
then all terms in the bracketed sum
are non-negative, because all the gamma functions
in fractions
$
\displaystyle{
\frac{\Gamma(a+1)}{\Gamma(b+1)\Gamma(a-b+1)}
=
\left(
\begin{array}{c}
a\\b
\end{array}
\right)
}
$
are non-negative,
i.e.,
$\delta_3\ge 0$.
If $\alpha_1<0$, then
$$
\delta_3\ge
-
\frac{\sin \pi\alpha_3}{\pi}
\sum\limits_{i_2=0}^{r_2} a_{i_2}(\lambda_2)
\frac{\Gamma(r_1+\alpha_1)}{\Gamma(r_1+1)|\Gamma(\alpha_1)|}
~
\frac{\Gamma(r_3+\alpha_3+1)\Gamma(n-i_2-r_1-r_3-\alpha_3)}{\Gamma(n-i_2-r_1+1)}
$$
$$
{}
\ge
-\frac{1}{3\pi r_1}
\sum\limits_{i_2=0}^{r_2}
\frac{\Gamma(r_2+\alpha_2+1)}{i_2!\Gamma(r_2+\alpha_2-i_2+1)}
~
%a_{i_2}(\lambda_2)
\frac{(r_3+1)!(r_2-i_2)!}{(r_2+r_3+1-i_2)!}
\ge
-\frac{1}{3\pi r_1}
\sum\limits_{i_2=0}^{r_2}
\left(
\begin{array}{c}
r_2+1\\
i_2
\end{array}
\right)
\ge
-\frac{2^{r_2-1}}{r_1}.
$$
\hfill$\Box$

\begin{lemen}
\label{delta_4}
%Имеет место неравенство
$\delta_4\ge-\frac{2^{r_2-1}}{r_1}$.
\end{lemen}

\textsl{Proof}.
Due to Lemma~\ref{vspom1} with
$m=n-i_2$, $q=r_3+1$, $r=r_1$, $x=r_1+\alpha_1$, $y=r_3+\alpha_3+1$
we obtain
$$
\delta_4=
\frac{|\sin\pi\alpha_1|\sin\pi\alpha_3}{\pi^2}
\sum\limits_{i_2=0}^{r_2-1}
a_{i_2}(\lambda_2)
\left(
\sum\limits_{i_3=r_3+1}^{n-r_1-1-i_2}
\frac{\Gamma(r_1+\alpha_1+1)\Gamma(i_1-r_1-\alpha_1)}{i_1!}
%\left|
%\frac{\Gamma(r_3+\alpha_3+1)}
%{i_3!~\Gamma(r_3+\alpha_3+1-i_3)}
%\right|
\frac{\Gamma(r_3+\alpha_3+1)\Gamma(i_3-r_3-\alpha_3)}{i_3!}
\right.
$$
$$
{}
-
\left.
\sum\limits_{i_3=r_3+2}^{n-r_1-i_2}
\frac{\Gamma(r_1+\alpha_1)\Gamma(i_1-r_1+1-\alpha_1)}{i_1!}
~
\frac{\Gamma(r_3+\alpha_3+2)\Gamma(i_3-r_3-1-\alpha_3)}{i_3!}
\right)
$$

$$
{}=
\frac{|\sin\pi\alpha_1|\sin\pi\alpha_3}{\pi^2}
\sum\limits_{i_2=0}^{r_2-1}
a_{i_2}(\lambda_2)
{\mathcal D}^{**}(n-i_2,r_3+1,r_1,r_1+\alpha_1,r_3+\alpha_3+1)
$$
$$
{}=
\frac{|\sin\pi\alpha_1|\sin\pi\alpha_3}{\pi^2}
\sum\limits_{i_2=0}^{r_2-1} a_{i_2}(\lambda_2)
\left(
\frac{1}
{
(r_1+\alpha_1)
\left(
\begin{array}{c}
r_1\\
r_1+\alpha_1
\end{array}
\right)
(r_3+\alpha_3+1)
\left(
\begin{array}{c}
n-i_2-r_1\\
r_3+\alpha_3+1
\end{array}
\right)
}
\right.
$$
$$
\left.
{}
-
\frac{1}
{
(r_1+\alpha_1)
\left(
\begin{array}{c}
n-i_2-r_3-1\\
r_1+\alpha_1
\end{array}
\right)
(r_3+\alpha_3+1)
\left(
\begin{array}{c}
r_3+1\\
r_3+\alpha_3+1
\end{array}
\right)
}\right)
$$

$$
{}\ge
-
\frac{|\sin\pi\alpha_1|\sin\pi\alpha_3}{\pi^2}
\sum\limits_{i_2=0}^{r_2-1}
\frac{a_{i_2}(\lambda_2)}
{
(r_1+\alpha_1)
\left(
\begin{array}{c}
n-i_2-r_3-1\\
r_1+\alpha_1
\end{array}
\right)
(r_3+\alpha_3+1)
\left(
\begin{array}{c}
r_3+1\\
r_3+\alpha_3+1
\end{array}
\right)
}
$$
\begin{equation}
\label{v1}
{}=
-
\frac{|\sin\pi\alpha_1|\sin\pi\alpha_3}{\pi^2}
\sum\limits_{i_2=0}^{r_2-1}
a_{i_2}(\lambda_2)
~
\frac{\Gamma(r_1+\alpha_1)\Gamma(n-i_2-r_3-r_1-\alpha_1)}
{\Gamma(n-i_2-r_3)}
~
\frac{\Gamma(r_3+\alpha_3+1)\Gamma(1-\alpha_3)}
{\Gamma(r_3+2)}~.
\end{equation}
Note that
\begin{equation}
\label{v2}
\frac{\sin\pi\alpha_3}{\pi}
~
\frac{\Gamma(r_3+\alpha_3+1)\Gamma(1-\alpha_3)}
{\Gamma(r_3+2)}
=
\frac{\Gamma(r_3+\alpha_3+1)}
{\Gamma(r_3+2)\Gamma(\alpha_3)}
\le
{\Gamma(r_3+2)}{\Gamma(r_3+2)\Gamma(1)}
\le 1.
\end{equation}
Consider $\mu$ defined by
$$
\mu
\overset{def}{=}
\frac{|\sin\pi\alpha_1|}{\pi}
~
\frac{\Gamma(r_1+\alpha_1)\Gamma(n-i_2-r_3-r_1-\alpha_1)}
{\Gamma(n-i_2-r_3)}
$$
and make transformations
$$
\mu
=
\frac{|\sin\pi\alpha_1|}{\pi}
~
\frac{\Gamma(r_1+\alpha_1)\Gamma(r_2+1-i_2-\alpha_1)}
{\Gamma(n-i_2-r_3)}
=
\frac{|\sin\pi\alpha_1|}{\pi}
~
\frac{\Gamma(\alpha_1)\Gamma(1-\alpha_1)}{(n-i_2-r_3-1)!}
\prod\limits_{s=0}^{r_1-1}(s+\alpha_1)
\prod\limits_{\kappa=1}^{r_2-i_2}(\kappa-\alpha_1)
$$
$$
{}=
\frac{(1+\alpha_1)|\alpha_1|(1-\alpha_1)
}{(n-i_2-r_3-1)!}
\prod\limits_{s=2}^{r_1-1}(s+\alpha_1)
\prod\limits_{\kappa=2}^{r_2-i_2}(\kappa-\alpha_1)
\le
\frac{1}{2(n-i_2-r_3-1)!}
\prod\limits_{s=2}^{r_1-1}(s+\alpha_1)
\prod\limits_{\kappa=2}^{r_2-i_2}(\kappa-\alpha_1).
$$
If $i_2=r_2-1$, then
\begin{equation}
\label{v3}
\mu
\le
\frac{(r_1+\alpha_1-1)\ldots (2+\alpha_1)}
{2(r_1+r_2-r_2+1)!}
\le
\frac{r_1(r_1-1)\ldots 3}
{2(r_1+1)!}
\le\frac{1}{4r_1}~.
\end{equation}
If $i_2\le r_2-2$, then
\begin{equation}
\label{v4}
\mu
\le
\frac{3~\ldots~(r_1+r_2-i_2-2)}
{2(r_1+r_2-i_2)!}
\le
\frac{1}{2r_1^2}
\le\frac{1}{4r_1} ~.
\end{equation}
It follows from~(\ref{v1}), (\ref{v2}), (\ref{v3}), (\ref{v4})
that
$$
\delta_4
\ge
-
\frac{1}{4r_1}
\sum\limits_{i_2=0}^{r_2} a_{i_2}(\lambda_2)
=
-
\frac{1}{4r_1}
\sum\limits_{i_2=0}^{r_2}
\left|
\frac{\Gamma(r_2+\alpha_2+1)}{i_2!~\Gamma(r_2+\alpha_2+1-i_2)}
\right|
\ge
-
\frac{1}{4r_1}
\sum\limits_{i_2=0}^{r_2}
\left(
\begin{array}{c}
r_2+1\\
i_2
\end{array}
\right)
\ge
-\frac{2^{r_2-1}}{r_1}.
$$
\hfill$\Box$

\begin{lemen}
\label{delta_6}
%Имеет место неравенство
$
\displaystyle{
%\begin{equation}
%\label{lemen_delta_6}
\delta_6
\ge
\sum\limits_{i_2=r_2+1}^{n-r_3-1}
a_{i_2}(\lambda_2)
\left(
\begin{array}{c}
r_1-1+\alpha_1\\
n-r_3-i_2-1
\end{array}
\right)
\left(
\begin{array}{c}
r_3+\alpha_3\\
r_3+1
\end{array}
\right).
%\end{equation}
}
$
\end{lemen}

\textsl{Proof}.
Using Lemma~\ref{vspom} with
$p=n-r_3-2-i_2$, $m=n-i_2$, $x=r_1+\alpha_1-1$,
$y=r_3+\alpha_3+1$
we obtain
$$
\delta_6=
\sum\limits_{i_2=r_2+1}^{n-r_3-1}
a_{i_2}(\lambda_2)
\sum\limits_{i_1=0}^{n-r_3-1-i_2}
\frac{\Gamma(r_1+\alpha_1+1)}{i_1!\left|\Gamma(r_1+\alpha_1-i_1+1)\right|}
~
\frac{\Gamma(r_3+\alpha_3+1)}{i_3!\left|\Gamma(r_3+\alpha_3-i_3+1)\right|}
$$
$$
{}-
\sum\limits_{i_2=r_2+1}^{n-r_3-2}
a_{i_2}(\lambda_2)
\sum\limits_{i_1=0}^{n-r_3-2-i_2}
\frac{\Gamma(r_1+\alpha_1)}{i_1!\left|\Gamma(r_1+\alpha_1-i_1)\right|}
~
\frac{\Gamma(r_3+\alpha_3+2)}{i_3!\left|\Gamma(r_3+\alpha_3-i_3+2)\right|}
$$
$$
{}=
\frac{\sin\pi\alpha_3}{\pi}
\sum\limits_{i_2=r_2+1}^{n-r_3-2}
a_{i_2}(\lambda_2)
\left(
\sum\limits_{i_1=0}^{n-r_3-1-i_2}
\frac{\Gamma(r_1+\alpha_1+1)}{i_1!\Gamma(r_1+\alpha_1-i_1+1)}
~
\frac{\Gamma(r_3+\alpha_3+1)\Gamma(i_3-r_3-\alpha_3)}{i_3!}
\right.
$$
$$
\left.
{}-
\sum\limits_{i_1=0}^{n-r_3-2-i_2}
\frac{\Gamma(r_1+\alpha_1)}{i_1!\Gamma(r_1+\alpha_1-i_1)}
~
\frac{\Gamma(r_3+\alpha_3+2)\Gamma(i_3-r_3-\alpha_3-1)}{i_3!}
\right)
$$
$$
{}+
\frac{\sin\pi\alpha_3}{\pi}
a_{n-r_3-1}(\lambda_2)
\frac{\Gamma(r_1+\alpha_1+1)}{0!~\Gamma(r_1+\alpha_1+1)}
~
\frac{\Gamma(r_3+\alpha_3+1)\Gamma(1-\alpha_3)}{(r_3+1)!}
$$

$$
{}=
\frac{\sin\pi\alpha_3}{\pi}
\sum\limits_{i_2=r_2+1}^{n-r_3-2}
a_{i_2}(\lambda_2)
{\mathcal D}^{*}(n-r_3-2-i_2, n-i_2, r_1+\alpha_1-1, r_3+\alpha_3+1)
$$
$$
{}+
\frac{\sin\pi\alpha_3}{\pi}
~
\frac{a_{n-r_3-1}(\lambda_2)
}
{
(r_3+1+\alpha_3)
\left(
\begin{array}{c}
r_3+1\\
r_3+1+\alpha_3
\end{array}
\right)
}
$$

$$
{}=
\frac{\sin\pi\alpha_3}{\pi}
\sum\limits_{i_2=r_2+1}^{n-r_3-2}
a_{i_2}(\lambda_2)
\left(
\sum\limits_{i_1=0}^{n-r_3-i_2-2}
\frac{2
\left(
\begin{array}{c}
r_1-1+\alpha_1\\
i_1
\end{array}
\right)
}
{(r_3+1+\alpha_3)
\left(
\begin{array}{c}
n-i_2-i_1\\
r_3+1+\alpha_3
\end{array}
\right)
}
\right.
$$
$$
{}
\left.
+
\frac{
\left(
\begin{array}{c}
r_1-1+\alpha_1\\
n-r_3-i_2-1
\end{array}
\right)
}
{(r_3+1+\alpha_3)
\left(
\begin{array}{c}
r_3+1\\
r_3+1+\alpha_3
\end{array}
\right)
}
\right)
+
\frac{\sin\pi\alpha_3}{\pi}
~
\frac{a_{n-r_3-1}(\lambda_2)
}
{
(r_3+1+\alpha_3)
\left(
\begin{array}{c}
r_3+1\\
r_3+1+\alpha_3
\end{array}
\right)
}
$$
$$
{}=
\frac{\sin\pi\alpha_3}{\pi}
\sum\limits_{i_2=r_2+1}^{n-r_3-2}
a_{i_2}(\lambda_2)
\sum\limits_{i_1=0}^{n-r_3-i_2-2}
\frac{2
\left(
\begin{array}{c}
r_1-1+\alpha_1\\
i_1
\end{array}
\right)
}
{(r_3+1+\alpha_3)
\left(
\begin{array}{c}
n-i_2-i_1\\
r_3+1+\alpha_3
\end{array}
\right)
}
$$
$$
{}
+
%\frac{\sin\pi\alpha_3}{\pi}
\sum\limits_{i_2=r_2+1}^{n-r_3-1}
a_{i_2}(\lambda_2)
\left(
\begin{array}{c}
r_1-1+\alpha_1\\
n-r_3-i_2-1
\end{array}
\right)
\left(
\begin{array}{c}
r_3+\alpha_3\\
r_3+1
\end{array}
\right).
$$
Since all the terms in the first sum are non-negative,
we obtain the assertion of the lemma.
\hfill$\Box$

\begin{lemen}
\label{delta_5}
$
\displaystyle{
\delta_5
\ge
-
\sum\limits_{i_2=r_2+1}^{n-r_1-1} a_{i_2}(\lambda_2)
\sum\limits_{i_3=0}^{n-r_1-1-i_2}
\left(
\begin{array}{c}
r_3+1\\
i_3
\end{array}
\right)
-
\sum\limits_{i_2=r_2+1}^{n-r_1} a_{i_2}(\lambda_2)
\left(
\begin{array}{c}
r_3+\alpha_3\\
n-r_1-i_2
\end{array}
\right)
\frac{\Gamma(r_1+\alpha_1)}{\Gamma(r_1+1)|\Gamma(\alpha_1)|}.
}
$
\end{lemen}

\textsl{Proof}.
Using Lemma~\ref{vspom} with
$p=n-r_1-i_2-1$, $m=n-i_2$, $x=r_3+\alpha_3$, $y=r_1+\alpha_1$
we obtain
$$
\delta_5=
\sum\limits_{i_2=r_2+1}^{n-r_1-1}
a_{i_2}(\lambda_2)
\sum\limits_{i_3=0}^{n-r_1-1-i_2}
\frac{\Gamma(r_1+\alpha_1+1)}{i_1!\left|\Gamma(r_1+\alpha_1-i_1+1)\right|}
~
\frac{\Gamma(r_3+\alpha_3+1)}{i_3!\left|\Gamma(r_3+\alpha_3-i_3+1)\right|}
$$
$$
{}-
\sum\limits_{i_2=r_2+1}^{n-r_1}
a_{i_2}(\lambda_2)
\sum\limits_{i_3=0}^{n-r_1-i_2}
\frac{\Gamma(r_1+\alpha_1)}{i_1!\left|\Gamma(r_1+\alpha_1-i_1)\right|}
~
\frac{\Gamma(r_3+\alpha_3+2)}{i_3!\left|\Gamma(r_3+\alpha_3-i_3+2)\right|}
%a_{i_1}(r_1-1,\alpha_1)a_{i_3}(r_3+1,\alpha_3)
$$

$$
{}=
-
\frac{|\sin\pi\alpha_1|}{\pi}
\sum\limits_{i_2=r_2+1}^{n-r_1-1}
a_{i_2}(\lambda_2)
\left(
\sum\limits_{i_3=0}^{n-r_1-i_2}
\frac{\Gamma(r_1+\alpha_1)\Gamma(i_1-r_1+1-\alpha_1)}{i_1!}
~
\frac{\Gamma(r_3+\alpha_3+2)}{i_3!~\Gamma(r_3+\alpha_3-i_3+2)}
\right.
$$
$$
\left.
-
\sum\limits_{i_3=0}^{n-r_1-1-i_2}
\frac{\Gamma(r_1+\alpha_1+1)\Gamma(i_1-r_1-\alpha_1)}{i_1!}
~
\frac{\Gamma(r_3+\alpha_3+1)}{i_3!~\Gamma(r_3+\alpha_3-i_3+1)}
\right)
$$
$$
{}-
\frac{|\sin\pi\alpha_1|}{\pi}
%\sum\limits_{i_2=n-r_1}^{n-r_1}
a_{n-r_1}(\lambda_2)
%\sum\limits_{i_3=0}^{n-r_1-i_2}
%a_{i_1}(\lambda_1)a_{i_3}(\lambda_3)
\frac{\Gamma(r_1+\alpha_1)\Gamma(1-\alpha_1)}{r_1!}
~
\frac{\Gamma(r_3+\alpha_3+2)}{0!~\Gamma(r_3+\alpha_3+2)}
$$

$$
{}=-
\frac{|\sin\pi\alpha_1|}{\pi}
\sum\limits_{i_2=r_2+1}^{n-r_1-1}
a_{i_2}(\lambda_2)
{\mathcal D}^{*}(n-r_1-i_2-1, n-i_2, r_3+\alpha_3, r_1+\alpha_1)
$$
$$
{}
-
\frac{|\sin\pi\alpha_1|}{\pi}~
a_{n-r_1}(\lambda_2)
\frac{\Gamma(r_1+\alpha_1)\Gamma(1-\alpha_1)}{r_1!}
$$

$$
{}=-
\frac{|\sin\pi\alpha_1|}{\pi}
\sum\limits_{i_2=r_2+1}^{n-r_1-1}
a_{i_2}(\lambda_2)
%\left(
\sum\limits_{i_3=0}^{n-r_1-i_2-1}
\frac{2
\left(
\begin{array}{c}
r_3+\alpha_3\\
i_3
\end{array}
\right)
}
{(r_1+\alpha_1)
\left(
\begin{array}{c}
n-i_2-i_3\\
r_1+\alpha_1
\end{array}
\right)
}
$$

$$
{}
-
\frac{|\sin\pi\alpha_1|}{\pi}
\sum\limits_{i_2=r_2+1}^{n-r_1}
a_{i_2}(\lambda_2)
\frac{
\left(
\begin{array}{c}
r_3+\alpha_3\\
n-r_1-i_2
\end{array}
\right)
}
{(r_1+\alpha_1)
\left(
\begin{array}{c}
r_1\\
r_1+\alpha_1
\end{array}
\right)
}.
$$
Now we can estimate $|\delta_5|$ as follows:
$$
|\delta_5|
\le
\frac{|\sin \pi\alpha_1|}{\pi}
\sum\limits_{i_2=r_2+1}^{n-r_1-1} a_{i_2}(\lambda_2)
\sum\limits_{i_3=0}^{n-r_1-1-i_2}
\frac{2
\left(
\begin{array}{c}
r_3+\alpha_3\\
i_3
\end{array}
\right)
}
{(r_1+\alpha_1)
\left(
\begin{array}{c}
r_1+1\\
r_1+\alpha_1
\end{array}
\right)
}
$$
$$
{}
+
\frac{|\sin\pi\alpha_1|}{\pi}~
\sum\limits_{i_2=r_2+1}^{n-r_1}
a_{i_2}(\lambda_2)
\left(
\begin{array}{c}
r_3+\alpha_3\\
n-r_1-i_2
\end{array}
\right)
\frac{\Gamma(r_1+\alpha_1)\Gamma(1-\alpha_1)}{r_1!}
$$

$$
{}=
\frac{|\sin \pi\alpha_1|}{\pi}
\sum\limits_{i_2=r_2+1}^{n-r_1-1} a_{i_2}(\lambda_2)
%\left(
\sum\limits_{i_3=0}^{n-r_1-1-i_2}
2
\left(
\begin{array}{c}
r_3+\alpha_3\\
i_3
\end{array}
\right)
\frac{\Gamma(r_1+\alpha_1)\Gamma(2-\alpha_1)}{\Gamma(r_1+2)}
%\right.
$$
$$
{}
+
%\frac{|\sin \pi\alpha_1|}{\pi}
\sum\limits_{i_2=r_2+1}^{n-r_1} a_{i_2}(\lambda_2)
\left(
\begin{array}{c}
r_3+\alpha_3\\
n-r_1-i_2
\end{array}
\right)
\frac{\Gamma(r_1+\alpha_1)}{\Gamma(r_1+1)|\Gamma(\alpha_1)|}
$$

$$
{}\le
%\frac{|\sin \pi\alpha_1|}{\pi}
\sum\limits_{i_2=r_2+1}^{n-r_1-1} a_{i_2}(\lambda_2)
\sum\limits_{i_3=0}^{n-r_1-1-i_2}
\left(
\begin{array}{c}
r_3+\alpha_3\\
i_3
\end{array}
\right)
\frac{2(1-\alpha_1)\Gamma(r_1+\alpha_1)}{(r_1+1)\Gamma(r_1+1)|\Gamma(\alpha_1)|}
$$
$$
{}
+
%\frac{|\sin \pi\alpha_1|}{\pi}
\sum\limits_{i_2=r_2+1}^{n-r_1} a_{i_2}(\lambda_2)
\left(
\begin{array}{c}
r_3+\alpha_3\\
n-r_1-i_2
\end{array}
\right)
\frac{\Gamma(r_1+\alpha_1)}{\Gamma(r_1+1)|\Gamma(\alpha_1)|}
$$

%$$
%{}\le
%\sum\limits_{i_2=r_2+1}^{n-r_1-1} a_{i_2}(\lambda_2)
%\sum\limits_{i_3=0}^{n-r_1-1-i_2}
%\left(
%\begin{array}{c}
%r_3+1\\
%i_3
%\end{array}
%\right)
%\frac{2(1-\alpha_1)
%%\Gamma(r_1+\alpha_1)
%}
%{(r_1+1)
%%\Gamma(r_1+1)
%|\Gamma(\alpha_1)|}
%+
%\sum\limits_{i_2=r_2+1}^{n-r_1} a_{i_2}(\lambda_2)
%\left(
%\begin{array}{c}
%r_3+\alpha_3\\
%n-r_1-i_2
%\end{array}
%\right)
%\frac{\Gamma(r_1+\alpha_1)}{\Gamma(r_1+1)|\Gamma(\alpha_1)|}
%$$

$$
{}\le
%\frac{1}{2}
\sum\limits_{i_2=r_2+1}^{n-r_1-1} a_{i_2}(\lambda_2)
\sum\limits_{i_3=0}^{n-r_1-1-i_2}
\left(
\begin{array}{c}
r_3+1\\
i_3
\end{array}
\right)
+
%\frac{1}{2}
\sum\limits_{i_2=r_2+1}^{n-r_1} a_{i_2}(\lambda_2)
\left(
\begin{array}{c}
r_3+\alpha_3\\
n-r_1-i_2
\end{array}
\right)
\frac{\Gamma(r_1+\alpha_1)}{\Gamma(r_1+1)|\Gamma(\alpha_1)|}.
$$
\hfill$\Box$

\begin{lemen}
\label{delta_2}
$$
\delta_2
=
-
\sum\limits_{i_2=r_2+1}^{n-r_3-1}
a_{i_2}(\lambda_2)
\left(
\begin{array}{c}
r_1+\alpha_1-1\\
n-r_3-1-i_2
\end{array}
\right)
\left(
\begin{array}{c}
r_3+\alpha_3\\
r_3+1
\end{array}
\right)
+
\sum\limits_{i_2=r_2+1}^{n-r_1}
a_{i_2}(\lambda_2)
%\sum\limits_{i_1=n-r_3-i_2+1}^{r_1}
\left(
\begin{array}{c}
r_1+\alpha_1-1\\
r_1
\end{array}
\right)
\left(
\begin{array}{c}
r_3+\alpha_3\\
n-r_1-i_2
\end{array}
\right).
$$
\end{lemen}

\textsl{Proof}.
Denote
$$
\delta_{2,k}={S}_{2,k}(r_1,r_2,r_3,\alpha_1,\alpha_2,\alpha_3)-{S}_{2,k}(r_1-1,r_2,r_3+1,\alpha_1,\alpha_2,\alpha_3),
\quad
1\le k\le 3.
$$
Using Lemma~\ref{vspom2} we obtain
$$
\delta_{2,1}=
\sum\limits_{i_2=n-r_1}^{n-r_3}
a_{i_2}(\lambda_2)
\sum\limits_{i_1=n-r_3-i_2}^{n-i_2}
\left(
\begin{array}{c}
r_1+\alpha_1\\
i_1
\end{array}
\right)
\left(
\begin{array}{c}
r_3+\alpha_3\\
i_3
\end{array}
\right)
$$
$$
{}
-
\sum\limits_{i_2=n-r_1+1}^{n-r_3-1}
a_{i_2}(\lambda_2)
\sum\limits_{i_1=n-r_3-1-i_2}^{n-i_2}
\left(
\begin{array}{c}
r_1+\alpha_1-1\\
i_1
\end{array}
\right)
\left(
\begin{array}{c}
r_3+\alpha_3+1\\
i_3
\end{array}
\right)
$$

$$
{}=
\sum\limits_{i_2=n-r_1+1}^{n-r_3-1}
a_{i_2}(\lambda_2)
\left(
\sum\limits_{i_1=n-r_3-i_2}^{n-i_2}
\left(
\begin{array}{c}
r_1+\alpha_1\\
i_1
\end{array}
\right)
\left(
\begin{array}{c}
r_3+\alpha_3\\
i_3
\end{array}
\right)
\right.
-
\left.
\sum\limits_{i_1=n-r_3-1-i_2}^{n-i_2}
\left(
\begin{array}{c}
r_1+\alpha_1-1\\
i_1
\end{array}
\right)
\left(
\begin{array}{c}
r_3+\alpha_3+1\\
i_3
\end{array}
\right)
\right)
$$
$$
{}+
a_{n-r_1}(\lambda_2)
\sum\limits_{i_1=r_1-r_3}^{r_1}
\left(
\begin{array}{c}
r_1+\alpha_1\\
i_1
\end{array}
\right)
\left(
\begin{array}{c}
r_3+\alpha_3\\
i_3
\end{array}
\right)
+
a_{n-r_3}(\lambda_2)
\sum\limits_{i_1=0}^{r_3}
\left(
\begin{array}{c}
r_1+\alpha_1\\
i_1
\end{array}
\right)
\left(
\begin{array}{c}
r_3+\alpha_3\\
i_3
\end{array}
\right)
$$

$$
{}=
\sum\limits_{i_2=n-r_1+1}^{n-r_3-1}
a_{i_2}(\lambda_2)
\left(
D^{***}(r_1+\alpha_1,r_3+\alpha_3,n-i_2,n-r_3-i_2,n-i_2)
-
\left(
\begin{array}{c}
r_1+\alpha_1-1\\
n-i_2
\end{array}
\right)
%\left( \begin{array}{c} r_3+\alpha_3+1\\ 0 \end{array} \right)
\right)
$$
$$
{}+
a_{n-r_1}(\lambda_2)
\sum\limits_{i_1=r_1-r_3}^{r_1}
\left(
\begin{array}{c}
r_1+\alpha_1\\
i_1
\end{array}
\right)
\left(
\begin{array}{c}
r_3+\alpha_3\\
i_3
\end{array}
\right)
+
a_{n-r_3}(\lambda_2)
\sum\limits_{i_1=0}^{r_3}
\left(
\begin{array}{c}
r_1+\alpha_1\\
i_1
\end{array}
\right)
\left(
\begin{array}{c}
r_3+\alpha_3\\
i_3
\end{array}
\right)
$$

$$
{}=
\sum\limits_{i_2=n-r_1+1}^{n-r_3-1}
a_{i_2}(\lambda_2)
\left(
\left(
\begin{array}{c}
r_1+\alpha_1-1\\
n-i_2
\end{array}
\right)
-
\left(
\begin{array}{c}
r_1+\alpha_1-1\\
n-r_3-i_2-1
\end{array}
\right)
\left(
\begin{array}{c}
r_3+\alpha_3\\
r_3+1
\end{array}
\right)
-
\left(
\begin{array}{c}
r_1+\alpha_1-1\\
n-i_2
\end{array}
\right)
%\left( \begin{array}{c} r_3+\alpha_3+1\\ 0 \end{array} \right)
\right)
$$
$$
{}+
a_{n-r_1}(\lambda_2)
\sum\limits_{i_1=r_1-r_3}^{r_1}
\left(
\begin{array}{c}
r_1+\alpha_1\\
i_1
\end{array}
\right)
\left(
\begin{array}{c}
r_3+\alpha_3\\
i_3
\end{array}
\right)
+
a_{n-r_3}(\lambda_2)
\sum\limits_{i_1=0}^{r_3}
\left(
\begin{array}{c}
r_1+\alpha_1\\
i_1
\end{array}
\right)
\left(
\begin{array}{c}
r_3+\alpha_3\\
i_3
\end{array}
\right)
$$

$$
{}=
-
\sum\limits_{i_2=n-r_1+1}^{n-r_3-1}
a_{i_2}(\lambda_2)
\left(
\begin{array}{c}
r_1+\alpha_1-1\\
n-r_3-i_2-1
\end{array}
\right)
\left(
\begin{array}{c}
r_3+\alpha_3\\
r_3+1
\end{array}
\right)
$$
$$
{}+
a_{n-r_1}(\lambda_2)
\sum\limits_{i_1=r_1-r_3}^{r_1}
\left(
\begin{array}{c}
r_1+\alpha_1\\
i_1
\end{array}
\right)
\left(
\begin{array}{c}
r_3+\alpha_3\\
i_3
\end{array}
\right)
+
a_{n-r_3}(\lambda_2)
\sum\limits_{i_1=0}^{r_3}
\left(
\begin{array}{c}
r_1+\alpha_1\\
i_1
\end{array}
\right)
\left(
\begin{array}{c}
r_3+\alpha_3\\
i_3
\end{array}
\right).
$$

Now we estumate $\delta_{2,2}$ and $\delta_{2,3}$:
$$
\delta_{2,2}=
\sum\limits_{i_2=n-r_3+1}^{n}
a_{i_2}(\lambda_2)
\sum\limits_{i_1=0}^{n-i_2}
\left(
\begin{array}{c}
r_1+\alpha_1\\
i_1
\end{array}
\right)
\left(
\begin{array}{c}
r_3+\alpha_3\\
i_3
\end{array}
\right)
-
\sum\limits_{i_2=n-r_3}^{n}
a_{i_2}(\lambda_2)
\sum\limits_{i_1=0}^{n-i_2}
\left(
\begin{array}{c}
r_1+\alpha_1-1\\
i_1
\end{array}
\right)
\left(
\begin{array}{c}
r_3+\alpha_3+1\\
i_3
\end{array}
\right)
$$

%$$
%{}=
%\sum\limits_{i_2=n-r_3+1}^{n}
%a_{i_2}(\lambda_2)
%\left(
%\sum\limits_{i_1=0}^{n-i_2}
%\left(
%\begin{array}{c}
%r_1+\alpha_1\\
%i_1
%\end{array}
%\right)
%\left(
%\begin{array}{c}
%r_3+\alpha_3\\
%i_3
%\end{array}
%\right)
%-
%%\sum\limits_{i_2=n-r_3}^{n}
%%a_{i_2}(\lambda_2)
%\sum\limits_{i_1=0}^{n-i_2}
%\left(
%\begin{array}{c}
%r_1+\alpha_1-1\\
%i_1
%\end{array}
%\right)
%\left(
%\begin{array}{c}
%r_3+\alpha_3+1\\
%i_3
%\end{array}
%\right)
%\right)
%$$
%$$
%{}
%-
%a_{n-r_3}(\lambda_2)
%\sum\limits_{i_1=0}^{r_3}
%\left(
%\begin{array}{c}
%r_1+\alpha_1-1\\
%i_1
%\end{array}
%\right)
%\left(
%\begin{array}{c}
%r_3+\alpha_3+1\\
%i_3
%\end{array}
%\right)
%$$

$$
{}=
\sum\limits_{i_2=n-r_3+1}^{n}
a_{i_2}(\lambda_2)
\left(
\sum\limits_{i_1=0}^{n-i_2}
\left(
\begin{array}{c}
r_1+\alpha_1-1\\
i_1
\end{array}
\right)
\left(
\begin{array}{c}
r_3+\alpha_3\\
i_3
\end{array}
\right)
+
\sum\limits_{i_1=1}^{n-i_2}
\left(
\begin{array}{c}
r_1+\alpha_1-1\\
i_1-1
\end{array}
\right)
\left(
\begin{array}{c}
r_3+\alpha_3\\
i_3
\end{array}
\right)
\right.
$$

$$
{}
\left.
-
\sum\limits_{i_1=0}^{n-i_2}
\left(
\begin{array}{c}
r_1+\alpha_1-1\\
i_1
\end{array}
\right)
\left(
\begin{array}{c}
r_3+\alpha_3\\
i_3
\end{array}
\right)
-
\sum\limits_{i_1=0}^{n-i_2-1}
\left(
\begin{array}{c}
r_1+\alpha_1-1\\
i_1
\end{array}
\right)
\left(
\begin{array}{c}
r_3+\alpha_3\\
i_3-1
\end{array}
\right)
\right)
$$
$$
{}
-
a_{n-r_3}(\lambda_2)
\sum\limits_{i_1=0}^{r_3}
\left(
\begin{array}{c}
r_1+\alpha_1-1\\
i_1
\end{array}
\right)
\left(
\begin{array}{c}
r_3+\alpha_3+1\\
i_3
\end{array}
\right)
=
-
a_{n-r_3}(\lambda_2)
\sum\limits_{i_1=0}^{r_3}
\left(
\begin{array}{c}
r_1+\alpha_1-1\\
i_1
\end{array}
\right)
\left(
\begin{array}{c}
r_3+\alpha_3+1\\
i_3
\end{array}
\right),
$$

$$
\delta_{2,3}=
\sum\limits_{i_2=r_2+1}^{n-r_1-1}
a_{i_2}(\lambda_2)
\sum\limits_{i_1=n-r_3-i_2}^{r_1}
\left(
\begin{array}{c}
r_1+\alpha_1\\
i_1
\end{array}
\right)
\left(
\begin{array}{c}
r_3+\alpha_3\\
i_3
\end{array}
\right)
$$
$$
{}
-
\sum\limits_{i_2=r_2+1}^{n-r_1}
a_{i_2}(\lambda_2)
\sum\limits_{i_1=n-r_3-i_2-1}^{r_1-1}
\left(
\begin{array}{c}
r_1+\alpha_1-1\\
i_1
\end{array}
\right)
\left(
\begin{array}{c}
r_3+\alpha_3+1\\
i_3
\end{array}
\right)
$$

$$
%\delta_{2,3}
{}=
\sum\limits_{i_2=r_2+1}^{n-r_1-1}
a_{i_2}(\lambda_2)
\left(
\sum\limits_{i_1=n-r_3-i_2}^{r_1}
\left(
\begin{array}{c}
r_1+\alpha_1\\
i_1
\end{array}
\right)
\left(
\begin{array}{c}
r_3+\alpha_3\\
i_3
\end{array}
\right)
-
\sum\limits_{i_1=n-r_3-i_2-1}^{r_1-1}
\left(
\begin{array}{c}
r_1+\alpha_1-1\\
i_1
\end{array}
\right)
\left(
\begin{array}{c}
r_3+\alpha_3+1\\
i_3
\end{array}
\right)
\right)
$$
$$
{}
-
a_{n-r_1}(\lambda_2)
\sum\limits_{i_1=r_1-r_3-1}^{r_1-1}
\left(
\begin{array}{c}
r_1+\alpha_1-1\\
i_1
\end{array}
\right)
\left(
\begin{array}{c}
r_3+\alpha_3+1\\
i_3
\end{array}
\right)
$$

$$
%\delta_{2,3}
{}=
\sum\limits_{i_2=r_2+1}^{n-r_1-1}
a_{i_2}(\lambda_2)
D^{***}(r_1+\alpha_1,r_3+\alpha_3,r_1,n-r_3-i_2,n-i_2)
$$
$$
{}
-
a_{n-r_1}(\lambda_2)
\sum\limits_{i_1=r_1-r_3-1}^{r_1-1}
\left(
\begin{array}{c}
r_1+\alpha_1-1\\
i_1
\end{array}
\right)
\left(
\begin{array}{c}
r_3+\alpha_3+1\\
i_3
\end{array}
\right)
$$

$$
{}=
\sum\limits_{i_2=r_2+1}^{n-r_1-1}
a_{i_2}(\lambda_2)
\left(
\left(
\begin{array}{c}
r_1+\alpha_1-1\\
r_1
\end{array}
\right)
\left(
\begin{array}{c}
r_3+\alpha_3\\
n-r_1-i_2
\end{array}
\right)
-
\left(
\begin{array}{c}
r_1+\alpha_1-1\\
n-r_3-i_2-1
\end{array}
\right)
\left(
\begin{array}{c}
r_3+\alpha_3\\
r_3+1
\end{array}
\right)
\right)
$$
$$
{}
-
a_{n-r_1}(\lambda_2)
\sum\limits_{i_1=r_1-r_3-1}^{r_1-1}
\left(
\begin{array}{c}
r_1+\alpha_1-1\\
i_1
\end{array}
\right)
\left(
\begin{array}{c}
r_3+\alpha_3+1\\
i_3
\end{array}
\right).
$$

Finally we have
$$
\delta_2=\delta_{2,1}+\delta_{2,2}+\delta_{2,3}
=
-
\sum\limits_{i_2=n-r_1+1}^{n-r_3-1}
a_{i_2}(\lambda_2)
\left(
\begin{array}{c}
r_1+\alpha_1-1\\
n-r_3-1-i_2
\end{array}
\right)
\left(
\begin{array}{c}
r_3+\alpha_3\\
r_3+1
\end{array}
\right)
$$
$$
{}+
a_{n-r_1}(\lambda_2)
\sum\limits_{i_1=r_1-r_3}^{r_1}
\left(
\begin{array}{c}
r_1+\alpha_1\\
i_1
\end{array}
\right)
\left(
\begin{array}{c}
r_3+\alpha_3\\
i_3
\end{array}
\right)
+
a_{n-r_3}(\lambda_2)
\sum\limits_{i_1=0}^{r_3}
\left(
\begin{array}{c}
r_1+\alpha_1\\
i_1
\end{array}
\right)
\left(
\begin{array}{c}
r_3+\alpha_3\\
i_3
\end{array}
\right)
$$
$$
{}
-
a_{n-r_3}(\lambda_2)
\sum\limits_{i_1=0}^{r_3}
\left(
\begin{array}{c}
r_1+\alpha_1-1\\
i_1
\end{array}
\right)
\left(
\begin{array}{c}
r_3+\alpha_3+1\\
i_3
\end{array}
\right)
$$
$$
{}+
\sum\limits_{i_2=r_2+1}^{n-r_1-1}
a_{i_2}(\lambda_2)
\left(
%\sum\limits_{i_1=n-r_3-i_2+1}^{r_1}
\left(
\begin{array}{c}
r_1+\alpha_1-1\\
r_1
\end{array}
\right)
\left(
\begin{array}{c}
r_3+\alpha_3\\
n-r_1-i_2
\end{array}
\right)
-
\left(
\begin{array}{c}
r_1+\alpha_1-1\\
n-r_3-i_2-1
\end{array}
\right)
\left(
\begin{array}{c}
r_3+\alpha_3\\
r_3+1
\end{array}
\right)
\right)
$$
$$
{}
-
a_{n-r_1}(\lambda_2)
\sum\limits_{i_1=r_1-r_3-1}^{r_1-1}
\left(
\begin{array}{c}
r_1+\alpha_1-1\\
i_1
\end{array}
\right)
\left(
\begin{array}{c}
r_3+\alpha_3+1\\
i_3
\end{array}
\right)
$$

$$
{}
=
-
\sum\limits_{i_2=n-r_1+1}^{n-r_3-1}
a_{i_2}(\lambda_2)
\left(
\begin{array}{c}
r_1+\alpha_1-1\\
n-r_3-1-i_2
\end{array}
\right)
\left(
\begin{array}{c}
r_3+\alpha_3\\
r_3+1
\end{array}
\right)
$$
$$
{}+
a_{n-r_1}(\lambda_2)
\left(
\sum\limits_{i_1=r_1-r_3}^{r_1}
\left(
\begin{array}{c}
r_1+\alpha_1\\
i_1
\end{array}
\right)
\left(
\begin{array}{c}
r_3+\alpha_3\\
r_1-i_1
\end{array}
\right)
-
%a_{n-r_1}(\lambda_2)
\sum\limits_{i_1=r_1-r_3-1}^{r_1-1}
\left(
\begin{array}{c}
r_1+\alpha_1-1\\
i_1
\end{array}
\right)
\left(
\begin{array}{c}
r_3+\alpha_3+1\\
r_1-i_1
\end{array}
\right)
\right)
$$
$$
{}
+
a_{n-r_3}(\lambda_2)
\left(
\sum\limits_{i_1=0}^{r_3}
\left(
\begin{array}{c}
r_1+\alpha_1\\
i_1
\end{array}
\right)
\left(
\begin{array}{c}
r_3+\alpha_3\\
r_3-i_1
\end{array}
\right)
-
%a_{n-r_3}(\lambda_2)
\sum\limits_{i_1=0}^{r_3}
\left(
\begin{array}{c}
r_1+\alpha_1-1\\
i_1
\end{array}
\right)
\left(
\begin{array}{c}
r_3+\alpha_3+1\\
r_3-i_1
\end{array}
\right)
\right)
$$
$$
{}+
\sum\limits_{i_2=r_2+1}^{n-r_1-1}
a_{i_2}(\lambda_2)
\left(
%\sum\limits_{i_1=n-r_3-i_2+1}^{r_1}
\left(
\begin{array}{c}
r_1+\alpha_1-1\\
r_1
\end{array}
\right)
\left(
\begin{array}{c}
r_3+\alpha_3\\
n-r_1-i_2
\end{array}
\right)
-
\left(
\begin{array}{c}
r_1+\alpha_1-1\\
n-r_3-i_2-1
\end{array}
\right)
\left(
\begin{array}{c}
r_3+\alpha_3\\
r_3+1
\end{array}
\right)
\right)
$$

%\sup\limits_{f\in C(\Delta) \atop {f\neq 0}}

$$
{}
=
-
\sum\limits_{{i_2=r_2+1} \atop i_2\neq n-r_1}^{n-r_3-1}
a_{i_2}(\lambda_2)
\left(
\begin{array}{c}
r_1+\alpha_1-1\\
n-r_3-1-i_2
\end{array}
\right)
\left(
\begin{array}{c}
r_3+\alpha_3\\
r_3+1
\end{array}
\right)
$$
$$
{}
+
a_{n-r_1}(\lambda_2)
D^{***}(r_1+\alpha_1,r_3+\alpha_3,r_1,r_1-r_3,r_1)
$$
$$
{}
+
a_{n-r_3}(\lambda_2)
\left(
%\left(
%\begin{array}{c}
%r_1+\alpha_1\\
%0
%\end{array}
%\right)
\left(
\begin{array}{c}
r_3+\alpha_3\\
r_3
\end{array}
\right)
-
\left(
\begin{array}{c}
r_1+\alpha_1-1\\
r_3
\end{array}
\right)
%\left(
%\begin{array}{c}
%r_3+\alpha_3+1\\
%0
%\end{array}
%\right)
+
D^{***}(r_1+\alpha_1,r_3+\alpha_3,r_3,1,r_3)
\right)
$$
$$
{}
+
\sum\limits_{i_2=r_2+1}^{n-r_1-1}
a_{i_2}(\lambda_2)
%\sum\limits_{i_1=n-r_3-i_2+1}^{r_1}
\left(
\begin{array}{c}
r_1+\alpha_1-1\\
r_1
\end{array}
\right)
\left(
\begin{array}{c}
r_3+\alpha_3\\
n-r_1-i_2
\end{array}
\right)
$$

%==================================================

$$
{}
=
-
\sum\limits_{{i_2=r_2+1} \atop i_2\neq n-r_1}^{n-r_3-1}
a_{i_2}(\lambda_2)
\left(
\begin{array}{c}
r_1+\alpha_1-1\\
n-r_3-1-i_2
\end{array}
\right)
\left(
\begin{array}{c}
r_3+\alpha_3\\
r_3+1
\end{array}
\right)
$$
$$
{}+
a_{n-r_1}(\lambda_2)
\left(
%\sum\limits_{i_1=r_1-r_3}^{r_1}
\left(
\begin{array}{c}
r_1+\alpha_1-1\\
r_1
\end{array}
\right)
\left(
\begin{array}{c}
r_3+\alpha_3\\
0
\end{array}
\right)
-
%a_{n-r_1}(\lambda_2)
%\sum\limits_{i_1=r_1-r_3-1}^{r_1-1}
\left(
\begin{array}{c}
r_1+\alpha_1-1\\
r_1-r_3-1
\end{array}
\right)
\left(
\begin{array}{c}
r_3+\alpha_3\\
r_3+1
\end{array}
\right)
\right)
$$
$$
{}
+
a_{n-r_3}(\lambda_2)
\left(
\left(
\begin{array}{c}
r_3+\alpha_3\\
r_3
\end{array}
\right)
-
\left(
\begin{array}{c}
r_1+\alpha_1-1\\
r_3
\end{array}
\right)
%\left(
%\begin{array}{c}
%r_3+\alpha_3+1\\
%0
%\end{array}
%\right)
+
%\sum\limits_{i_1=0}^{r_3}
\left(
\begin{array}{c}
r_1+\alpha_1-1\\
r_3
\end{array}
\right)
%\left( \begin{array}{c} r_3+\alpha_3\\0 \end{array} \right)
-
%\left( \begin{array}{c} r_1+\alpha_1-1\\ 0 \end{array} \right)
\left(
\begin{array}{c}
r_3+\alpha_3\\
r_3
\end{array}
\right)
\right)
$$
$$
{}+
\sum\limits_{i_2=r_2+1}^{n-r_1-1}
a_{i_2}(\lambda_2)
%\sum\limits_{i_1=n-r_3-i_2+1}^{r_1}
\left(
\begin{array}{c}
r_1+\alpha_1-1\\
r_1
\end{array}
\right)
\left(
\begin{array}{c}
r_3+\alpha_3\\
n-r_1-i_2
\end{array}
\right)
$$

%==================================================

$$
{}
=
-
\sum\limits_{i_2=r_2+1}^{n-r_3-1}
a_{i_2}(\lambda_2)
\left(
\begin{array}{c}
r_1+\alpha_1-1\\
n-r_3-1-i_2
\end{array}
\right)
\left(
\begin{array}{c}
r_3+\alpha_3\\
r_3+1
\end{array}
\right)
+
\sum\limits_{i_2=r_2+1}^{n-r_1}
a_{i_2}(\lambda_2)
%\sum\limits_{i_1=n-r_3-i_2+1}^{r_1}
\left(
\begin{array}{c}
r_1+\alpha_1-1\\
r_1
\end{array}
\right)
\left(
\begin{array}{c}
r_3+\alpha_3\\
n-r_1-i_2
\end{array}
\right).
$$
\hfill$\Box$

\section{Reducing to the case $r_1=n-1-r_2$, $r_3=0$}

\begin{lemen}
\label{shift_r3}
Let
$r_s\in {\Bbb Z}_+$, $1\le s\le 3$,
$r_1+r_2+r_3=n-1$,
$-1 < \alpha_1< 1$, $0\le \alpha_2,\alpha_3<1$,
$r_1-1+\alpha_1\ge r_2+\alpha_2$,
$r_1-1+\alpha_1\ge r_3+1+\alpha_3$.
Then
\begin{equation}
\label{lemen_shift_r3}
\mathfrak{L}_n(r_1-1,r_2,r_3+1,\alpha_1,\alpha_2,\alpha_3)
{}
\le
\mathfrak{L}_n(r_1,r_2,r_3,\alpha_1,\alpha_2,\alpha_3)
+
2^{r_2+r_3+2}
+\frac{2^{r_2}}{r_1}
-1
+
\left\{
\begin{array}{lll}
2^{r_3} & \mbox{if} & r_2\ge 1,\\
2^{r_3+1}\ln n & \mbox{if} & r_2=0.
\end{array}
\right.
\end{equation}
\end{lemen}

\textsl{Proof}.
Recall that $r_1\ge r_3+2$, $n\ge 4$.
Consider the difference
$$
\omega(r_1,r_2,r_3)
\overset{def}{=}
\mathfrak{L}_n(r_1,r_2,r_3,\alpha_1,\alpha_2,\alpha_3)
-
\mathfrak{L}_n(r_1-1,r_2,r_3+1,\alpha_1,\alpha_2,\alpha_3).
$$
Due to~(\ref{Lambda_S}), (\ref{delta_k})
and Lemmas~\ref{delta_1}, \ref{delta_3}, \ref{delta_4}, \ref{delta_6}, \ref{delta_5}, \ref{delta_2},
we obtain
$$
\omega(r_1,r_2,r_3)
=
\sum\limits_{s=1}^6
\delta_s
\ge
-2^{r_2+r_3+2}
+2^{r_3+1}
-\frac{2^{r_2-1}}{r_1}
-\frac{2^{r_2-1}}{r_1}
+
\sum\limits_{i_2=r_2+1}^{n-r_3-1}
a_{i_2}(\lambda_2)
\left(
\begin{array}{c}
r_1-1+\alpha_1\\
n-r_3-i_2-1
\end{array}
\right)
\left(
\begin{array}{c}
r_3+\alpha_3\\
r_3+1
\end{array}
\right)
%-
%\frac{(r_3+1)2^{r_3}}{\pi}
$$
$$
{}
-
\sum\limits_{i_2=r_2+1}^{n-r_1-1} a_{i_2}(\lambda_2)
\sum\limits_{i_3=0}^{n-r_1-1-i_2}
\left(
\begin{array}{c}
r_3+1\\
i_3
\end{array}
\right)
-
%\frac{1}{2}
\sum\limits_{i_2=r_2+1}^{n-r_1} a_{i_2}(\lambda_2)
\left(
\begin{array}{c}
r_3+\alpha_3\\
n-r_1-i_2
\end{array}
\right)
\frac{\Gamma(r_1+\alpha_1)}{\Gamma(r_1+1)|\Gamma(\alpha_1)|}
$$
$$
-
\sum\limits_{i_2=r_2+1}^{n-r_3-1}
a_{i_2}(\lambda_2)
\left(
\begin{array}{c}
r_1+\alpha_1-1\\
n-r_3-1-i_2
\end{array}
\right)
\left(
\begin{array}{c}
r_3+\alpha_3\\
r_3+1
\end{array}
\right)
+
\sum\limits_{i_2=r_2+1}^{n-r_1}
a_{i_2}(\lambda_2)
%\sum\limits_{i_1=n-r_3-i_2+1}^{r_1}
\left(
\begin{array}{c}
r_1+\alpha_1-1\\
r_1
\end{array}
\right)
\left(
\begin{array}{c}
r_3+\alpha_3\\
n-r_1-i_2
\end{array}
\right)
$$

$$
{}=
-2^{r_2+r_3+2}
+2^{r_3+1}
%-\frac{1}{r_1}
-\frac{2^{r_2}}{r_1}
-
\sum\limits_{i_2=r_2+1}^{n-r_1-1} a_{i_2}(\lambda_2)
\sum\limits_{i_3=0}^{n-r_1-1-i_2}
\left(
\begin{array}{c}
r_3+1\\
i_3
\end{array}
\right)
$$
$$
{}
-
%\frac{1}{2}
\sum\limits_{i_2=r_2+1}^{n-r_1} a_{i_2}(\lambda_2)
\left(
\begin{array}{c}
r_3+\alpha_3\\
n-r_1-i_2
\end{array}
\right)
\frac{\Gamma(r_1+\alpha_1)}{\Gamma(r_1+1)}
\left(
\frac{1}{|\Gamma(\alpha_1)|}
-
\frac{1}{\Gamma(\alpha_1)}
\right).
$$
If $0<\alpha_1<1$, then
$
\displaystyle{
\frac{1}{|\Gamma(\alpha_1)|}
-
\frac{1}{\Gamma(\alpha_1)}=0
}$.
If $-1<\alpha_1<0$, then
$
\displaystyle{
\frac{1}{|\Gamma(\alpha_1)|}
-
\frac{1}{\Gamma(\alpha_1)}
\le
1}$.
Thus
$$
\omega(r_1,r_2,r_3)
\ge
-2^{r_2+r_3+2}
+2^{r_3+1}
%-\frac{1}{r_1}
-\frac{2^{r_2}}{r_1}
-
\sum\limits_{i_2=r_2+1}^{n-r_1-1} a_{i_2}(\lambda_2)
\sum\limits_{i_3=0}^{n-r_1-1-i_2}
\left(
\begin{array}{c}
r_3+1\\
i_3
\end{array}
\right)
-
%\frac{1}{2}
\sum\limits_{i_2=r_2+1}^{n-r_1} a_{i_2}(\lambda_2)
\left(
\begin{array}{c}
r_3+\alpha_3\\
n-r_1-i_2
\end{array}
\right)
$$
$$
{}
\ge
-2^{r_2+r_3+2}
+2^{r_3+1}
%-\frac{1}{r_1}
-\frac{2^{r_2}}{r_1}
-
\sum\limits_{i_2=r_2+1}^{n-r_1} a_{i_2}(\lambda_2)
%\left(
\sum\limits_{i_3=0}^{n-r_1-i_2}
\left(
\begin{array}{c}
r_3+1\\
i_3
\end{array}
\right)
$$
%$$
%{}
%\ge
%%\frac{1}{2}
%-2^{r_2+r_3+2}
%%-\frac{1}{r_1}
%-\frac{2^{r_2+1}}{r_1}
%-
%2^{r_3+1}
%%\sum\limits_{i_3=0}^{r_3}
%%\left(
%%\begin{array}{c}
%%r_3+1\\
%%i_3
%%\end{array}
%%\right)
%\sum\limits_{i_2=r_2+1}^{n-r_1}
%\left|
%\frac{\Gamma(r_2+\alpha_2+1)}
%     {i_2!\Gamma(r_2+\alpha_2-i_2+1)}
%\right|
%$$
$$
{}\ge
-2^{r_2+r_3+2}
+2^{r_3+1}
-\frac{2^{r_2}}{r_1}
-
(2^{r_3+1}-1)~
\frac{\sin\pi\alpha_2}{\pi}
%\sum\limits_{i_3=0}^{r_3}
%\left(
%\begin{array}{c}
%r_3+1\\
%i_3
%\end{array}
%\right)
\sum\limits_{i_2=r_2+1}^{n-r_1}
\frac{\Gamma(r_2+\alpha_2+1)\Gamma(i_2-r_2-\alpha_2)}
     {i_2!}~.
$$
If $i_2=r_2+1$, then
$$
\frac{\sin\pi\alpha_2}{\pi}
~
\frac{\Gamma(r_2+\alpha_2+1)\Gamma(i_2-r_2-\alpha_2)}{i_2!}
=
\frac{\Gamma(r_2+\alpha_2+1)}{(r_2+1)!~\Gamma(\alpha_2)}
\le
\frac{\Gamma(r_2+2)}{(r_2+1)!~\Gamma(1)}
\le 1.
$$
Let $r_2\neq 0$.
If $i_2\ge r_2+2$,
then taking into account the fact that
$(1+\alpha_2)\alpha_2(1-\alpha_2)\le 0.5$
for $0\le \alpha_2\le 1$, we obtain
$$
\frac{\sin\pi\alpha_2}{\pi}
~
\Gamma(r_2+\alpha_2+1)\Gamma(i_2-r_2-\alpha_2)
=
\frac{\sin\pi\alpha_2}{\pi}
\Gamma(\alpha_2)\Gamma(1-\alpha_2)
\prod\limits_{s=0}^{r_2}
(s+\alpha_2)
\prod\limits_{\kappa=1}^{i_2-r_2-1}
(\kappa-\alpha_2)
$$
$$
{}=
(1+\alpha_2)\alpha_2(1-\alpha_2)
\prod\limits_{s=2}^{r_2}
(s+\alpha_2)
\prod\limits_{\kappa=2}^{i_2-r_2-1}
(\kappa-\alpha_2)
\le
\frac{2}{2}
\prod\limits_{s=2}^{r_2}
(s+\alpha_2)
\prod\limits_{\kappa=3}^{i_2-r_2-1}
(\kappa-\alpha_2)
\le
\prod\limits_{s=3}^{r_2+1}
s
\prod\limits_{\kappa=3}^{i_2-r_2-1}
\kappa
\le
\prod\limits_{\kappa=3}^{i_2-2}
\kappa.
$$
Then
$$
\frac{\sin\pi\alpha_2}{\pi}
\sum\limits_{i_2=r_2+1}^{n-r_1}
\frac{\Gamma(r_2+\alpha_2+1)\Gamma(i_2-r_2-\alpha_2)}{i_2!}
\le
\sum\limits_{i_2=r_2+1}^{n-r_1}
\frac{1}{i_2!}
\prod\limits_{\kappa=3}^{i_2-2}
\kappa
=
\sum\limits_{i_2=r_2+1}^{n-r_1}
\frac{1}{2(i_2-1)i_2}
$$
$$
{}
\le
1+
\frac{1}{2}\sum\limits_{i_2=r_2+2}^{n-r_1}
\frac{1}{(i_2-1)i_2}
=
1
+
\frac{1}{2(r_2+1)}
-
\frac{1}{2(n-r_1)}
\le
1
+
\frac{1}{2(r_2+1)}
\le\frac{3}{2}~,
$$
and
$$
\omega(r_1,r_2,r_3)
\ge
-2^{r_2+r_3+2}
+2^{r_3+1}
-\frac{2^{r-2}}{r_1}
-(2^{r_3+1}-1)\frac{3}{2}
\ge
-2^{r_2+r_3+2}
-\frac{2^{r-2}}{r_1}
-2^{r_3}+1.
$$
If $r_2=0$, then
$$
\frac{\sin\pi\alpha_2}{\pi}
\sum\limits_{i_2=r_2+1}^{n-r_1}
\frac{\Gamma(r_2+\alpha_2+1)\Gamma(i_2-r_2-\alpha_2)}{i_2!}
=
\frac{\sin\pi\alpha_2}{\pi}
\sum\limits_{i_2=1}^{n-r_1}
\frac{\Gamma(1+\alpha_2)\Gamma(i_2-\alpha_2)}{i_2!}
$$
$$
{}
\le
1+
\frac{\sin\pi\alpha_2}{\pi}
\sum\limits_{i_2=2}^{n-r_1}
\frac{\alpha_2\Gamma(\alpha_2)\Gamma(1-\alpha_2)(1-\alpha_2)\ldots(i_2-1-\alpha_2)}{i_2!}
$$
$$
\le
1+
\sum\limits_{i_2=2}^{n-r_1}
\frac{\alpha_2(1-\alpha_2)\ldots(i_2-1-\alpha_2)}{i_2!}
\le
1+
\frac{1}{2}
\sum\limits_{i_2=2}^{n-r_1}
\frac{1}{i_2}
\le 1+ \frac{1}{4}\ln (n-r_1)
\le
\ln n,
%1+\frac{\ln n}{4}~.
$$
and
$$
\omega(r_1,0,r_3)
\ge
-2^{r_2+r_3+2}
+2^{r_3+1}
-\frac{2^{r-2}}{r_1}
-(2^{r_3+1}-1)\ln 2
\ge
-2^{r_2+r_3+2}
-\frac{2^{r-2}}{r_1}
-2^{r_3+1}\ln 2+1.
$$

Finally we have
$$
\omega(r_1,r_2,r_3)
\ge
-2^{r_2+r_3+2}
%-\frac{1}{r_1}
-\frac{2^{r_2}}{r_1}
+1
-
\left\{
\begin{array}{lll}
2^{r_3} & \mbox{if} & r_2\ge 1,\\
2^{r_3+1}\ln n & \mbox{if} & r_2=0,
\end{array}
\right.
$$
and it follows from this that~(\ref{lemen_shift_r3}) holds.
\hfill$\Box$

\begin{lemen}
\label{zero_r3}
Let
$r_s\in {\Bbb Z}_+$, $1\le s\le 3$,
$r_1+r_2+r_3=n-1$,
$-1 < \alpha_1< 1$, $0\le \alpha_2,\alpha_3<1$,
$r_1+\alpha_1\ge r_2+\alpha_2$,
$r_1+\alpha_1\ge r_3+\alpha_3$.
Then
$$
\mathfrak{L}_n(r_1,r_2,r_3,\alpha_1,\alpha_2,\alpha_3)
\le
\mathfrak{L}_n(r_1+r_3,r_2,0,\alpha_1,\alpha_2,\alpha_3)
+
2^{r_2+r_3+2}
+
2^{r_2}
-r_3
+
\left\{
\begin{array}{lll}
2^{r_3} & \mbox{if} & r_2\ge 1,\\
2^{r_3+1}\ln n & \mbox{if} & r_2=0.
\end{array}
\right.
$$
\end{lemen}

\textsl{Proof}.
Let $r_2\ge 1.$
Applying $r_3$ times lemma~\ref{shift_r3}, we obtain
$$
\mathfrak{L}_n(r_1,r_2,r_3,\alpha_1,\alpha_2,\alpha_3)
\le
\mathfrak{L}_n(r_1+r_3,r_2,0,\alpha_1,\alpha_2,\alpha_3)
+
\sum\limits_{s=0}^{r_3-1}
2^{r_2+s+2}
+
\sum\limits_{s=0}^{r_3-1}
\frac{2^{r_2}}{r_1}
+
\sum\limits_{s=0}^{r_3-1} 2^s
-\sum\limits_{s=0}^{r_3-1}1
$$
$$
{}\le
\mathfrak{L}_n(r_1+r_3,r_2,0,\alpha_1,\alpha_2,\alpha_3)
%{\mathcal{L}}_n(r_1+r_3+\alpha_1,r_2+\alpha_2,\alpha_3)
+
2^{r_2+r_3+2}
%+ 2^{r_2+1}
+
2^{r_2}
+
2^{r_3}-r_3.
$$
The case $r_2=0$ is treated similarly.
\hfill$\Box$

\section{Reducing to the case $r_1=n-1$, $r_2=r_3=0$}

\begin{lemen}
\label{zero_r2_r3}
Let
$r_s\in {\Bbb Z}_+$, $1\le s\le 3$,
$r_3=0$,
$r_1+r_2=n-1$,
$-1 < \alpha_1< 1$, $0\le \alpha_2,\alpha_3<1$,
$r_1+\alpha_1\ge r_2+\alpha_2$.
Then
$$
\mathfrak{L}_n(r_1,r_2,0,\alpha_1,\alpha_2,\alpha_3)
\le
\mathfrak{L}_n(r_1+r_2,0,0,\alpha_1,\alpha_2,\alpha_3)
+
2^{r_2+2}
+1
-r_2
+
2^{r_2+1}\ln n.
$$
\end{lemen}

\textsl{Proof}.
We treat this case similarly to Section~5,
where the parameter $r_2$ plays the role of the parameter $r_3$.

\hfill$\Box$

\begin{teoen}
\label{zero_zero}
Let
$r_s\in {\Bbb Z}_+$, $1\le s\le 3$,
$r_1+r_2+r_3=n-1$,
$-1 < \alpha_1< 1$, $0\le \alpha_2,\alpha_3<1$,
$r_1+\alpha_1\ge r_2+\alpha_2$,
$r_1+\alpha_1\ge r_3+\alpha_3$.
Then
$$
\mathfrak{L}_n(r_1,r_2,r_3,\alpha_1,\alpha_2,\alpha_3)
\le
\mathfrak{L}_n(n-1,0,0,\alpha_1,\alpha_2,\alpha_3)
+
2^{2n/3}(10+2\ln n).
$$
\end{teoen}

\textsl{Proof}.
We apply Lemmas~\ref{zero_r3} and~\ref{zero_r2_r3}.
If $r_2\neq 0$, then
$$
\mathfrak{L}_n(r_1,r_2,r_3,\alpha_1,\alpha_2,\alpha_3)
-
\mathfrak{L}_n(r_1+r_2+r_3,0,0,\alpha_1,\alpha_2,\alpha_3)
\le
2^{r_2+r_3+2}
+2^{r_2}
+2^{r_3}
+2^{r_2+2}
+1
+2^{r_2+1}\ln n
-(r_2+r_3)
$$
$$
{}\le
4\cdot 2^{r_2+r_3}
+2^{r_2}(1+4+2\ln n)
+2^{r_3}
\le
4\cdot 2^{r_2+r_3}+2^{r_2+r_3}(6+2\ln n)
\le
2^{2n/3}(10+2\ln n).
$$
If $r_2=0$, then
$$
\mathfrak{L}_n(r_1,r_2,r_3,\alpha_1,\alpha_2,\alpha_3)
-
\mathfrak{L}_n(r_1+r_2+r_3,0,0,\alpha_1,\alpha_2,\alpha_3)
\le
2^{r_3+2}+1+2^{r_3+1}\ln n
\le
2^{2n/3}(10+2\ln n).
$$
\hfill$\Box$

\section{The upper bound of $\Lambda_n$}

{\bf 8.1. Statement of the main result.}
The purpose of this section is to prove the following main theorem of this article.

\begin{teoen}
\label{teo_osn}
Let $n\ge 4$. Then
\begin{equation}\label{teo_teo_osn}
\Lambda_n
=
\Lambda_{n,d}
\le
(7+\mu_n)
\frac{2^{n+1}}{en(\ln n -\ln 2)}
\left(1+\frac{15}{n-3}\right),
\end{equation}
where
$$
\mu_n
\le
\frac{3en(\ln n)^3}{2^{n/3}}
+
\frac{en^2\ln n}{2^n}
\to 0\quad
\mbox{as}
\quad
n\to \infty.
$$
\end{teoen}

Let $k,m,p\in{\Bbb Z}_+$. If $k,p$ are odd,
we assume that
$$
\sum\limits_{s=k/2}^n
=
\sum\limits_{s=(k+1)/2}^n
\quad
\mbox{and}
\quad
\sum\limits_{s=m}^{p/2}
=
\sum\limits_{s=m}^{(p-1)/2}.
$$

\textsl{Proof}.
Lemma~\ref{restrict} and Theorem~\ref{zero_zero} show that
\begin{equation}
\label{ddd}
\Lambda_n
\le
\underset
{
-1\le\alpha_1\le1
\atop
{
0\le \alpha_2,\alpha_3\le 1
\atop
0\le\alpha_1+\alpha_2+\alpha_3\le1
}
}
\max
\mathfrak{L}_n(n-1,0,0,\alpha_1,\alpha_2,\alpha_3)
+
2^{2n/3}(10+2\ln n).
\end{equation}

To estimate $\mathfrak{L}_n(n-1,0,0,\alpha_1,\alpha_2,\alpha_3)$
let's estimate each term in~(\ref{Lambda_S}) for $r_1=n-1$, $r_2=r_3=0$.
Denote
$\bar{{\mathcal S}}_k=S_k(n-1,0,0,\alpha_1,\alpha_2,\alpha_3)$,
$1\le k\le 6$. Then
\begin{equation}
\label{vvv}
\mathfrak{L}_n(n-1,0,0,\alpha_1,\alpha_2,\alpha_3)
=
\sum\limits_{k=1}^6
\bar{{\mathcal S}}_k.
\end{equation}

Everywhere we assume that
$i=(i_1,i_2,i_3)\in I$, $i_1,i_2,i_3\in {\Bbb Z}_+$, $i_1+i_2+i_3=n,$
$\lambda_1=(n-1+\alpha_1)/n$, $\lambda_2=\alpha_2/n$, $\lambda_3=\alpha_3/n$,
$\lambda=(\lambda_1,\lambda_2,\lambda_3)$,
$-1<\alpha_1<1$, $0<\alpha_2,\alpha_3<1$, $\alpha_1+\alpha_2+\alpha_3=1$.
Now we will deal with the sums $\bar{{\mathcal S}}_k$.

{\bf 8.2. The upper bound of $\bar{{\mathcal S}}_6=S_6(n-1,0,0,\alpha_1,\alpha_2,\alpha_3)$.}
Denote
$$
\varepsilon(i_1)
=\varepsilon
=
\left\{
\begin{array}{lll}
0,&\mbox{if}& \frac{n-i_1}{2}\in\Bbb{N},\\
1,&\mbox{if}& \frac{n-i_1}{2}\notin\Bbb{N}.
\end{array}
\right.
$$
If $\frac{n-i_1}{2}\in\Bbb{Z}_+$ then the polynomial
$l_i(\lambda)$ such that 
$
\displaystyle{
i=
\left(i_1,\frac{n-i_1}{2},\frac{n-i_1}{2}
\right)
}$ denote $\mathfrak{p} (i_1)$, i.e.,
$$
\mathfrak{p}(i_1)
=
l_i(\lambda)\Big|_{i=\left(i_1,\frac{n-i_1}{2},\frac{n-i_1}{2}\right)}.
$$

Consider $\bar{{\mathcal S}}_6$
and represent it as a sum of terms
$$
\bar{{\mathcal S}}_6
=
\sum\limits_{s=1}^4
\sigma_s,
$$
where
$$
\sigma_1
=
\sum\limits_{i_1=0}^{n-4}
\left(
\sum\limits_{i_2=1}^{\frac{n-i_1-1}{2}}
|l_i(\lambda)|
+\frac{1-\varepsilon(i_1)}{2}\mathfrak{p}(i_1)
\right),
\quad
\sigma_2
=
\sum\limits_{i_1=0}^{n-4}
\left(
\sum\limits_{i_3=1}^{\frac{n-i_1-1}{2}}
|l_i(\lambda)|
+\frac{1-\varepsilon(i_1)}{2}\mathfrak{p}(i_1)
\right),
$$
%$$
%\sigma_3
%=
%\sum\limits_{i_1=n-3}^{n-2}
%\sum\limits_{i_2=\frac{n-i_1-1}{2}}^{\frac{n-i_1+1}{2}}
%%\sum\limits_{\frac{n-i_1-2}{2}<i_2<\frac{n-i_1+2}{2}}
%|l_i(\lambda)|.
%%|l_i(\lambda)|\Big|_{i=(n-2,1,1)}.
%$$
$$
\sigma_3
=
|l_i(\lambda)|\Big|_{i=(n-3,1,2)}
+
|l_i(\lambda)|\Big|_{i=(n-3,2,1)},
\quad
\sigma_4
=
|l_i(\lambda)|\Big|_{i=(n-2,1,1)}.
$$

Before estimating the summands of $\sigma_s$,
we prove several auxiliary lemmas.
The following auxiliary lemma is similar
to the lemma from~\cite{Ture1940},
the proof is given here for the convenience of the reader.

\begin{lemen}
\label{Ture_1}
Let
$$
\Phi(\alpha,m)
=
\frac{\sin \pi\alpha}{\pi}
~
\frac{\Gamma(1+\alpha)\Gamma(m-\alpha)}{m!},
$$
where $0<\alpha<1$, $m$ is an integer, $m\ge 2$.
Then
$$
\Phi(\alpha,m)
\le\frac{\alpha
(1-\alpha)
2^\alpha}{m^{1+\alpha}}.
$$
If $m=1$, then $\Phi(\alpha,m)=\Phi(\alpha,1)=\alpha$.
\end{lemen}

\textsl{Proof}~\cite{Ture1968}.
For $m=1,2$ the inequality can be checked directly by substituting $m$.
Let $m\ge 3$.
First, we estimate the value
$$
\varphi(\alpha,m)
\overset{def}{=}
\prod\limits_{s=2}^{m-1}(s-\alpha).
$$
Take the logarithm and make the following transformations:
$$
\ln \varphi(\alpha,m)
=
\sum\limits_{s=2}^{m-1}\ln(s-\alpha)
=
\sum\limits_{s=2}^{m-1}\ln s
-
\sum\limits_{s=2}^{m-1}
(\ln s-\ln(s-\alpha))
=
\ln (m-1)!
-\sum\limits_{s=2}^{m-1}\ln\left(1+\frac{\alpha}{s-\alpha}\right)
$$
$$
{}
\le
\ln (m-1)!
-
\sum\limits_{s=2}^{m-1}\frac{1}{1+\frac{\alpha}{s-\alpha}}~\frac{\alpha}{s-\alpha}
=
\ln (m-1)!
-
\alpha\sum\limits_{s=2}^{m-1}\frac{1}{s}
\le
\ln (m-1)!
-
\alpha\int\limits_{s=2}^{m}\frac{ds}{s}
=
\ln (m-1)!-\alpha\ln m+\alpha\ln 2.
$$
So,
$
\displaystyle{
\varphi(\alpha,m)
\le
\frac{(m-1)!2^\alpha}{m^\alpha}.
}
$
Then
$$
\Phi(\alpha,m)
=
\frac{\alpha\Gamma(m-\alpha)}{m!~\Gamma(1-\alpha)}
=
\frac{\alpha}{m!}
\prod\limits_{s=1}^{m-1}(s-\alpha)
=
\frac{\alpha(1-\alpha)}{m!}\varphi(\alpha,m)
\le
\frac{\alpha(1-\alpha)(m-1)!2^\alpha}{m!~m^\alpha}
=
\frac{\alpha(1-\alpha)2^\alpha}{m^{1+\alpha}}~.
$$
\hfill$\Box$

\begin{lemen}
\label{ln}
Let $\alpha\ge 1$, $m> 1$,
$
\displaystyle{
\phi(\alpha,m)
=
\frac{\alpha}{m^{\alpha}}.
}
$
Then
$
\displaystyle{
\phi(\alpha,m)\le\frac{1}{e\ln m}~.
}
$
\end{lemen}

\textsl{Proof}.
Consider the function $\phi(\alpha,m)$ as a function
of $\alpha$ on $(0,+\infty)$.
Since
$
\phi'_\alpha=
m^{-\alpha}(1-\alpha\ln m)
$,
the point $\alpha=1/\ln m$ is the maximum point of the function
$\phi$.
Then
$$
\phi(\alpha,m)
\le
\phi\left(\frac{1}{\ln m}~,m\right)
=
\frac{1}{\ln m ~ m^{1/\ln m}}
=\frac{1}{e\ln m}~.
$$
\hfill$\Box$

%================================================

\begin{lemen}
\label{sum}
Let $n\ge 4$. Then
\begin{equation}
\label{lemen_sum}
\sum\limits_{s=1}^{n}
\frac{2^s}{s}
\le
\frac{2^{n+1}}{n}+\frac{2^{n+1}}{(n-3)^2}~.
\end{equation}
\end{lemen}

\textsl{Proof}.
For $n=\overline{4,8}$ the inequality~(\ref{lemen_sum})
is verified by direct calculations.
Let $n\ge 9$. We proceed by induction.
Suppose the inequality~(\ref{lemen_sum})
to be true for $n=\mathfrak{t}-1$, i.e.,
$$
\sum\limits_{s=1}^{\mathfrak{t}-1}
\frac{2^s}{s}
\le
\frac{2^{\mathfrak{t}}}{\mathfrak{t}-1}+\frac{2^{\mathfrak{t}}}{(\mathfrak{t}-4)^2}
\overset{def}{=}\varphi_1(\mathfrak{t}).
$$
Now prove that~(\ref{lemen_sum}) holds for $n=\mathfrak{t}$.
This is equivalent to the inequality
$$
\sum\limits_{s=1}^{\mathfrak{t}-1}
\frac{2^s}{s}
\le
\frac{2^{\mathfrak{t}}}{\mathfrak{t}}+\frac{2^{\mathfrak{t}+1}}{(\mathfrak{t}-3)^2}
\overset{def}{=}\varphi_2(\mathfrak{t}).
$$
Thus, it suffices for us to prove that
$\varphi_1(\mathfrak{t})\le\varphi_2(\mathfrak{t})$.
Consider the difference
$$
\varphi_2(\mathfrak{t})-\varphi_1(\mathfrak{t})
=
\frac{2^\mathfrak{t}(3\mathfrak{t}^3-40\mathfrak{t}^2+145\mathfrak{t}-144)}
{\mathfrak{t}(\mathfrak{t}-1)(\mathfrak{t}-3)^2(\mathfrak{t}-4)^2}
=
\frac{2^\mathfrak{t}g(\mathfrak{t})}{\mathfrak{t}(\mathfrak{t}-1)(\mathfrak{t}-3)^2(\mathfrak{t}-4)^2}~,
$$
where
$g(\mathfrak{t})=3\mathfrak{t}^3-40\mathfrak{t}^2+145\mathfrak{t}-144$.
Since $g'(\mathfrak{t})>0$ for $\mathfrak{t}\ge 7$ and $g(9)>0$,
then $g(\mathfrak{t})>0$ for all $\mathfrak{t}\ge 9$,
which implies that $\varphi_1(\mathfrak{t})\le\varphi_2(\mathfrak{t})$
for all $\mathfrak{t}\ge 9$.
\hfill$\Box$

Denote
$$
\nu_p(s)=
\left\{
\begin{array}{lll}
\alpha_p &\mbox{if} & s=1,\\
\alpha_p(1-\alpha_p)2^{\alpha_p} & \mbox{if} & s\ge 2,
\end{array}
\right.
$$
where $p=2,3$.

\begin{lemen}
\label{fm}
Let $i_2\ge 1$, $i_3\ge 1$, $i_2+i_3\ge 3$. Then
\begin{equation}\label{fm1}
|l_i(\lambda)|
\le
\frac{(1-\alpha_3)
%\alpha_2 (1-\alpha_2) 2^{\alpha_2}
}{e(\ln n-\ln 2)}
\frac{\Gamma(n+1)}{i_1!~\Gamma(n-i_1)}
~\frac{
\nu_2(i_2)
}
{i_2^{1+\alpha_2}i_3(i_3-1)}
\quad (i_2\ge 1, i_3\ge 2),
%\frac{\Gamma(i_3-1)}{i_3!}~,
\end{equation}
\begin{equation}\label{fm2}
|l_i(\lambda)|
\le
\frac{(1-\alpha_2)
%\alpha_2 (1-\alpha_2) 2^{\alpha_2}
}{e(\ln n-\ln 2)}\frac{\Gamma(n+1)}{i_1!~\Gamma(n-i_1)}
~\frac{
\nu_3(i_3)
}{i_3^{1+\alpha_3}i_2(i_2-1)}
\quad(i_2\ge 2, i_3\ge 1).
%\frac{\Gamma(i_2-1)}{i_2!}~.
\end{equation}
\end{lemen}

\textsl{Proof}.
We will prove~(\ref{fm1}).
Let $i_2,i_3\ge 1$.
Due to~(\ref{abs_l}) and
%presentation of functions $|l_i(\lambda)|$,
Lemmas~\ref{monot}, \ref{Ture_1}
we have
$$
l_i(\lambda)
=
\prod\limits_{s=1}^3 a_{i_s}(\lambda_s)
=
\frac{\sin \pi\alpha_2\sin\pi\alpha_3}{\pi^2}
\frac{\Gamma(n+\alpha_1)}{i_1!~\Gamma(n+\alpha_1-i_1)}
\frac{\Gamma(1+\alpha_2)\Gamma(i_2-\alpha_2)}{i_2!}
\frac{\Gamma(1+\alpha_3)\Gamma(i_3-\alpha_3)}{i_3!}
$$

$$
{}=
\frac{\sin\pi\alpha_3}{\pi}
%\sum\limits_{i_1=0}^{n-4}
\frac{\Gamma(n+1-\alpha_2-\alpha_3)}{i_1!~\Gamma(n+1-i_1-\alpha_2-\alpha_3)}
%\sum\limits_{i_2=1}^{(n-i_1)/2-1}
\frac{\sin \pi\alpha_2}{\pi}
\frac{\Gamma(1+\alpha_2)\Gamma(i_2-\alpha_2)}{i_2!}
%\cdot
\frac{\Gamma(1+\alpha_3)\Gamma(i_3-\alpha_3)}{i_3!}
$$
$$
{}\le
\frac{\sin\pi\alpha_3}{\pi}
%\sum\limits_{i_1=0}^{n-4}
\frac{\Gamma(n+1-\alpha_3)}{i_1!~\Gamma(n+1-i_1-\alpha_3)}
%\sum\limits_{i_2=1}^{(n-i_1)/2-1}
~\frac{
%\alpha_2 (1-\alpha_2) 2^{\alpha_2}
\nu_2(i_2)
}{i_2^{1+\alpha_2}}~
%\cdot
\frac{\Gamma(1+\alpha_3)\Gamma(i_3-\alpha_3)}{i_3!}
$$
$$
{}=
\frac{\sin\pi\alpha_3}{\pi}
\frac{\Gamma(1+\alpha_3)\Gamma(n+1-\alpha_3)}{(n+1)!}
\frac{(n+1)!}{i_1!}
~\frac{
nu_2(i_2)
}{i_2^{1+\alpha_2}}~
\frac{
%\Gamma(n+1)
\Gamma(i_3-\alpha_3)}
{
%i_1!~
i_3!~\Gamma(n+1-i_1-\alpha_3)}
$$
$$
{}
\le
\frac{\alpha_3(1-\alpha_3)2^{\alpha_3}}{(n+1)^{1+\alpha_3}}
\frac{(n+1)!}{i_1!}
~\frac{
\nu_2(i_2)
}{i_2^{1+\alpha_2}}~
\frac{
%\Gamma(n+1)
\Gamma(i_3-1)}
{
%i_1!~
i_3!~\Gamma(n+1-i_1-1)}.
$$
By Lemma~\ref{ln}, we have
$$
\frac{\alpha_3 2^{\alpha_3}}{(n+1)^{\alpha_3}}
=
\phi\left(\alpha_3,\frac{n+1}{2}\right)
\le\frac{1}{e(\ln (n+1) -\ln 2)}
\le
\frac{1}{e(\ln n -\ln 2)}
%\le\frac{1}{e\ln n}\left(1+\frac{1}{\ln n-1}\right).
$$
which implies
\begin{equation}
\label{i_3_ge_1}
|l_i(\lambda)|
\le
\frac{1-\alpha_3}{e(\ln n-\ln 2)}
\frac{\Gamma(n+1)}{i_1!}
~\frac{
\nu_2(i_2)
}{i_2^{1+\alpha_2}}~.
\frac{
\Gamma(i_3-\alpha_3)}
{
i_3!~\Gamma(n+1-i_1-\alpha_3)}.
\end{equation}

If $i_3\ge 2$, $i_2\ge 1$, then due to~(\ref{i_3_ge_1}) and Lemma~\ref{monot},
we obtain
$$
l_i(\lambda)
\le
\frac{1-\alpha_3}{e(\ln n-\ln 2)}
\frac{\Gamma(n+1)}{i_1!}
~\frac{
\nu_2(i_2)
}{i_2^{1+\alpha_2}}~
\frac{\Gamma(i_3-1)}{i_3!~\Gamma(n-i_1)}
=
\frac{1-\alpha_3}{e(\ln n-\ln 2)}
\frac{\Gamma(n+1)}{i_1!~\Gamma(n-i_1)}
~\frac{
\nu_2(i_2)
}{i_2^{1+\alpha_2}i_3(i_3-1)}~.
$$
The case of~(\ref{fm2}) is similar.
\hfill$\Box$

\begin{lemen}
\label{pol_p}
Let $(n-i_1)$ be even, $i_1\le n-4$.
Then
\begin{equation}
\label{lemen_pol_p}
\mathfrak{p}(i_1)
\le
\frac{1}{e(\ln n-1)}
\frac{\Gamma(n+1)}{i_1!~\Gamma(n-i_1)}
\frac{4\alpha_2(1-\alpha_2)}{(n-i_1)(n-i_1-1)}.
\end{equation}
\end{lemen}

\textsl{Proof}.
Consider $\mathfrak{p}(i_1)$ and apply Lemma~\ref{fm}.
We also take into account that
$i_1\le n-4$. Then
$$
\mathfrak{p}(i_1)
\le
\frac{1
}{e(\ln n-1)}
\frac{\Gamma(n+1)}{i_1!~\Gamma(n-i_1)}
~\frac{
\alpha_2(1-\alpha_2)2^{\alpha_2}
(1-\alpha_3)
}{i_2^{1+\alpha_2}i_3(i_3-1)}
\Big|_{i_2=i_3=\frac{n-i_1}{2}}
$$
$$
{}
=
\frac{1
}{e(\ln n-1)}
\frac{\Gamma(n+1)}{i_1!~\Gamma(n-i_1)}
\frac{
8\alpha_2(1-\alpha_2)
(1-\alpha_3)4^{\alpha_2}}
{(n-i_1)^{2+\alpha_2}(n-i_1-2)}
$$

$$
{}
\le
\frac{1}{e(\ln n-1)}
%\left(
%\begin{array}{c}
%n\\
%n-i_1
%\end{array}
%\right)
\frac{\Gamma(n+1)}{i_1!~\Gamma(n-i_1)}
\frac{1}{(n-i_1)^2}
\frac{
8\alpha_2(1-\alpha_2)(1-\alpha_3)4^{\alpha_2}}
{4^{\alpha_2}2}
$$

$$
{}
\le
\frac{1}{e(\ln n-1)}
\frac{\Gamma(n+1)}{i_1!~\Gamma(n-i_1)}
\frac{4\alpha_2(1-\alpha_2)}{(n-i_1)(n-i_1-1)}.
$$
\hfill$\Box$

%=======================================================

\begin{lemen}
\label{sigma_12}
The following inequality holds:
\begin{equation}
\label{lemen_sigma_12}
\sigma_{s}
\le
\frac{C_0
}{e(\ln n-1)}
\sum\limits_{i_1=0}^{n-4}
\left(
\begin{array}{c}
n\\
n-i_1
\end{array}
\right)
\frac{1}{n-i_1-1},
\quad s=1,2,
\end{equation}
where
$$
C_0=\frac{(\ln 2)^2+12\ln 2+28}{4\ln 2+12}<2.5.
$$
\end{lemen}

\textsl{Proof}.
We will prove~(\ref{lemen_sigma_12}) for $s=1$.
Let $\tilde{\mathfrak{p}}(i_1)$ be a function such that
$$
\mathfrak{p}(i_1)
=
\frac{1}{e(\ln n-1)}
\frac{n!}{i_1!~(n-1-i_1)!}
\tilde{\mathfrak{p}}(i_1).
$$
Using Lemma~\ref{fm} we have
\begin{equation}
\label{raz}
\sigma_{1}
\le
\frac{1
%(1-\alpha_3)
}{e(\ln n-1)}
\sum\limits_{i_1=0}^{n-4}
%\frac{\Gamma(n+1)}{i_1!~\Gamma(n-i_1)}
\frac{n!}{i_1!~(n-1-i_1)!}
\left(
\frac{\alpha_2}{(n-i_1-1)(n-i_1-2)}
\right.
$$
$$
{}
\left.
+
\sum\limits_{i_2=2}^{(n-i_1-1)/2}
\frac{\alpha_2 (1-\alpha_2) 2^{\alpha_2}}{i_2^{1+\alpha_2}i_3(i_3-1)}
+
\frac{1-\varepsilon(i_1)}{2}\tilde{\mathfrak{p}}(i_1)
\right).
\end{equation}
Now we estimate the terms in brackets as follows:
$$
h_1(i_1)
\overset{def}{=}
\frac{\alpha_2}{(n-i_1-1)(n-i_1-2)}
=
\frac{\alpha_2}{(n-i_1)}
\left(
\frac{1}{n-i_1-1}+\frac{2}{(n-i_1-1)(n-i_1-2)}
\right)
$$
%Если $i_1\le n-4$, то
\begin{equation}
\label{h1}
%h_1(i_1)
\le
\frac{\alpha_2}{(n-i_1)}
\left(
\frac{1}{n-i_1-1}+\frac{2}{2(n-i_1-1)}
\right)
=
\frac{2\alpha_2}{(n-i_1)(n-i_1-1)},
\end{equation}
%(при написании неравенства учитываем тот факт, что
%$i_1\le n-3$, если $n-i_1$ нечетное, и $i_i\le n-4$,
%если $n-i_1$ четное),
%Если $i_1=n-3$, то
%\begin{equation}
%\label{h1_dop}
%h_1(n-3)
%=
%\frac{3\alpha_2}{(n-i_1)(n-i_1-1)}
%\Big|_{i_1=n-3}.
%\end{equation}
%Второе слагаемое - это

$$
h_2(i_1)
\overset{def}{=}
\sum\limits_{i_2=2}^{(n-i_1-1)/2}
\frac{\alpha_2 (1-\alpha_2) 2^{\alpha_2}}{i_2^{1+\alpha_2}i_3(i_3-1)}
%$$
%Так как при $i_1=n-3$ верхний предел суммирования меньше нижнего
%предела, то $h_2(n-3)=0$. Оценим $h_2(i_1)$ при $i_1\le n-4$
%следующим образом:
%$$
=
\sum\limits_{i_2=2}^{(n-i_1-1)/2}
\frac{\alpha_2 (1-\alpha_2) 2^{\alpha_2}}{i_2^{1+\alpha_2}(n-i_1-i_2)(n-i_1-i_2-1)}
$$
$$
{}=
\sum\limits_{i_2=2}^{(n-i_1-1)/2}
\frac{\alpha_2 (1-\alpha_2) 2^{\alpha_2}}{n-i_1}
\left(
\frac{1+(n-i_1-i_2)-(n-i_1-i_2)}{i_2^{\alpha_2}(n-i_1-i_2)(n-i_1-i_2-1)}
+\frac{1}{i_2^{1+\alpha_2}(n-i_1-i_2-1)}
\right)
$$
$$
{}\le
\alpha_2 (1-\alpha_2) 2^{\alpha_2}
\sum\limits_{i_2=2}^{(n-i_1-1)/2}
\left(
\frac{1}{n-i_1}
\left(
\frac{1}{2^{\alpha_2}(n-i_1-i_2-1)}
-\frac{1}{2^{\alpha_2}(n-i_1-i_2)}
\right)
\right.
$$
$$
{}+
%\sum\limits_{i_2=2}^{(n-i_2-2)/2}
\left.
\frac{1}{(n-i_1)(n-i_1-1)}
\left(
\frac{1}{i_2^{\alpha_2}(n-i_1-i_2-1)}
+\frac{1}{i_2^{1+\alpha_2}}
\right)
\right)
$$

$$
{}\le
\frac{\alpha_2 (1-\alpha_2) 2^{\alpha_2}}{n-i_1}
\int\limits_{2}^{(n-i_1+\varepsilon)/2}
\left(
\frac{1}{2^{\alpha_2}(n-i_1-t-1)}
-\frac{1}{2^{\alpha_2}(n-i_1-t)}
\right)dt
$$

$$
{}+
%\sum\limits_{i_2=2}^{(n-i_2-2)/2}
\frac{\alpha_2 (1-\alpha_2) 2^{\alpha_2}}{(n-i_1)(n-i_1-1)}
\left(
\int\limits_{2}^{(n-i_1+\varepsilon)/2}
\frac{dt}{2^{\alpha_2}(n-i_1-t-1)}
+
\int\limits_{1}^{(n-i_1-2+\varepsilon)/2}
\frac{dt}{t^{1+\alpha_2}}
\right)
$$

$$
{}=
\frac{\alpha_2 (1-\alpha_2) 2^{\alpha_2}}{2^{\alpha_2}(n-i_1)}
\left(
-\ln\frac{n-i_1-2-\varepsilon}{2}+\ln(n-i_1-3)
+\ln\frac{n-i_1-\varepsilon}{2}-\ln(n-i_1-2)
\right)
$$
$$
{}+
\frac{\alpha_2 (1-\alpha_2) 2^{\alpha_2}}{2^{\alpha_2}(n-i_1)(n-i_1-1)}
\left(
-\ln\frac{n-i_1-2-\varepsilon}{2}+\ln(n-i_1-3)
%+\ln\frac{n-i_1-2}{2}%-\ln(n-i_1-1)
-\frac{2^{\alpha_2}~2^{\alpha_2}}{\alpha_2(n-i_1-2+\varepsilon)^{\alpha_2}}
+\frac{1}{\alpha_2}
\right)
$$

$$
{}\le
\frac{\alpha_2 (1-\alpha_2) }{n-i_1}
\left(
-\ln(n-i_1-2-\varepsilon)+\ln(n-i_1-3)
+\ln(n-i_1-\varepsilon)-\ln(n-i_1-2)
\right)
$$
$$
{}
+
\frac{\alpha_2 (1-\alpha_2)}{(n-i_1)(n-i_1-1)}
\left(
\ln 2
+\frac{1}{\alpha_2}
%-\frac{2^{\alpha_2}\cdot2^{\alpha_2}}{\alpha_2(n-i_1-2)^{\alpha_2}}
\right).
$$
If $\varepsilon=0$, then
$$
-\ln(n-i_1-2-\varepsilon)+\ln(n-i_1-3)
+\ln(n-i_1-\varepsilon)-\ln(n-i_1-2)
$$
$$
{}=
-\ln(n-i_1-2)^2+\ln(n-i_1-3)
+\ln(n-i_1)
$$
$$
{}
\le
-\ln(n-i_1-1)(n-i_1-3)
+\ln(n-i_1-3)
+\ln(n-i_1)
=
-\ln(n-i_1-1)
+\ln(n-i_1)
$$
$$
{}
=
\ln\left(1+\frac{1}{n-i_1-1}\right)
\le
\frac{1}{n-i_1-1}.
$$
If $\varepsilon=1$, then
$$
-\ln(n-i_1-2-\varepsilon)+\ln(n-i_1-3)
+\ln(n-i_1-\varepsilon)-\ln(n-i_1-2)
$$
$$
{}=
-\ln(n-i_1-3)+\ln(n-i_1-3)
+\ln(n-i_1-1)
-\ln(n-i_1-2)
=
\ln\left(1+\frac{1}{n-i_1-2}\right)
$$
$$
{}
\le
\frac{1}{n-i_1-2}
=
\frac{1}{n-i_1-1}
\left(
1+
\frac{1}{n-i_1-2}
\right)
\le
\frac{3}{2(n-i_1-1)}.
$$
Thus
\begin{equation}
\label{h2}
h_2(i_1)
\le
\frac{\alpha_2 (1-\alpha_2)}{(n-i_1)(n-i_1-1)}
\left(
1+\frac{\varepsilon(i_1)}{2}
+
\ln 2
+\frac{1}{\alpha_2}
%-\frac{2^{\alpha_2}\cdot2^{\alpha_2}}{\alpha_2(n-i_1-2)^{\alpha_2}}
\right).
\end{equation}
Further, taking into account~(\ref{lemen_pol_p}), we have
\begin{equation}
\label{h3}
h_3(i_1)
\overset{def}{=}
\frac{1-\varepsilon(i_1)}{2}\tilde{\mathfrak{p}}(i_1)
\le
\frac{1-\varepsilon(i_1)}{2}
\frac{4\alpha_2(1-\alpha_2)}{(n-i_1)(n-i_1-1)}.
\end{equation}

It follows from~(\ref{raz}), (\ref{h1}), (\ref{h2}), (\ref{h3})
that
$$
\sigma_{1}
\le
\frac{
}{e(\ln n-1)}
\sum\limits_{i_1=0}^{n-4}
\frac{n!}{i_1!~(n-1-i_1)!}
\frac{1}{(n-i_1)(n-i_1-1)}
$$

$$
{}\times
\left(
2\alpha_2
+
\alpha_2 (1-\alpha_2)
\left(
1+\frac{\varepsilon(i_1)}{2}
+
\ln 2
+\frac{1}{\alpha_2}
\right)
+
4\alpha_2(1-\alpha_2)
\frac{1-\varepsilon(i_1)}{2}
\right)
$$

$$
{}
=
\frac{1
}{e(\ln n-1)}
\sum\limits_{i_1=0}^{n-4}
\left(
\begin{array}{c}
n\\
n-i_1
\end{array}
\right)
\frac{1}{n-i_1-1}
\left(
2\alpha_2
+
\alpha_2 (1-\alpha_2)
\left(
3+\ln 2
-\frac{3\varepsilon(i_1)}{2}
+\frac{1}{\alpha_2}
\right)
\right)
$$

\begin{equation}\label{sig1}
{}
\le
\frac{1
}{e(\ln n-1)}
\sum\limits_{i_1=0}^{n-4}
\left(
\begin{array}{c}
n\\
n-i_1
\end{array}
\right)
\frac{2\alpha_2+(1-\alpha_2)
\left( (3+\ln 2)\alpha_2+1 \right)
}
{n-i_1-1}.
\end{equation}
Consider the function
$$
\eta(\alpha_2)=2\alpha_2+(1-\alpha_2)
\left( (3+\ln 2)\alpha_2+1 \right).
$$
The point
$
\displaystyle{
\alpha_2=\frac{4+\ln 2}{6+2\ln 2}
}$
is the maximum point of this function, which implies that
\begin{equation}
\label{sig2}
\eta(\alpha_2)
\le
\eta
\left(
\frac{4+\ln 2}{6+2\ln 2}
\right)
=
\frac{(\ln 2)^2+12\ln 2+28}{4\ln 2+12}<2.5.
\end{equation}
From~(\ref{sig1}) and ~(\ref{sig2}) we obtain~(\ref{lemen_sigma_12}) for $s=1$.
The case $s=2$ is proved in the same way.
\hfill$\Box$

\begin{lemen}
\label{sigma_3}
$
\displaystyle{
\sigma_3\le
\frac{3}{e(\ln n-1)}
\left(
\begin{array}{c}
n\\
n-3
\end{array}
\right)
\frac{1}{2}.
%\end{equation}
}
$
\end{lemen}

\textsl{Proof}.
Consider $\sigma_3$ and make estimates using Lemma~\ref{fm}:
$$
\sigma_3\le
\frac{1-\alpha_3}{e(\ln n-1)}
\frac{\Gamma(n+1)}{(n-3)!~\Gamma(3)}
\frac{\alpha_2}{2\cdot 1}
+
\frac{1-\alpha_2}{e(\ln n-1)}
\frac{\Gamma(n+1)}{(n-3)!~\Gamma(3)}
\frac{\alpha_3}{2\cdot 1}
$$
\begin{equation}
\label{ccc}
{}=
\frac{1}{e(\ln n-1)}
\frac{\Gamma(n+1)}{(n-3)!~\Gamma(4)}~
\frac{3}{2}~
(\alpha_2(1-\alpha_3)+\alpha_3(1-\alpha_2)).
\end{equation}
The maximum of the function
$\zeta(\alpha_2,\alpha_3)
=
\alpha_2(1-\alpha_3)+\alpha_3(1-\alpha_2)$,
$0\le\alpha_2,\alpha_3\le 1$,
equals $1$.
Then from~(\ref{ccc})
we obtain the assertion of the lemma.
\hfill$\Box$

\begin{lemen}
\label{sigma_4}
$
\displaystyle{
\sigma_4\le
\frac{2}{e(\ln n-1)}
\left(
\begin{array}{c}
n\\
n-2
\end{array}
\right).
}
$
\end{lemen}

\textsl{Proof}.
Due to the estimate~(\ref{i_3_ge_1})
found in the proof of the lemma~\ref{fm} we obtain
$$
\sigma_4
=
|l_i(\lambda)|\Big|_{i=(n-2,1,1)}
\le
\frac{1-\alpha_3}{e(\ln n-1)}
\frac{\Gamma(n+1)}{(n-2)!}
~\frac{
\alpha_2
}{1^{1+\alpha_2}}~.
\frac{
\Gamma(1-\alpha_3)}
{
1!~\Gamma(n+1-(n-2)-\alpha_3)}
$$
$$
{}
=
\frac{\alpha_2(1-\alpha_3)}{e(\ln n-1)}
\frac{n!}{(n-2)!}~\frac{2}{2}
\le
\frac{2}{e(\ln n-1)}
\left(
\begin{array}{c}
n\\
n-2
\end{array}
\right).
$$
\hfill$\Box$

To estimate the sum $\bar{\mathcal{S}}_6$
we prove the following lemma.
A similar result was proved in~\cite{Ture1968}.
The proposed proof is given
for the convenience of the reader.

\begin{lemen}
\label{Ture_2}
Let $n\ge 4$. Then
\begin{equation}
\label{lemen_Ture_2}
\sum\limits_{i_1=0}^{n-2}
\frac{n!}{i_1!~(n-i_1)!}
~
\frac{1}{n-i_1-1}
=
\sum\limits_{s=2}^{n}
\left(
\begin{array}{c}
n\\
s
\end{array}
\right)
\frac{1}{s-1}
\le
\frac{2^{n+1}}{n}\left(1+\frac{15}{n-3}\right).
\end{equation}
\end{lemen}

\textsl{Proof}~\cite{Ture1968}.
First note that
$$
\Upsilon(n)
\overset{def}{=}
\sum\limits_{i_1=0}^{n-2}
\frac{n!}{i_1!~(n-i_1)!}
~
\frac{1}{n-i_1-1}
=
\sum\limits_{s=2}^{n}
\left(
\begin{array}{c}
n\\
s
\end{array}
\right)
\frac{1}{s-1}
=
\sum\limits_{s=2}^{n}
\int\limits_0^1
\left(
\begin{array}{c}
n\\
s
\end{array}
\right)
t^{s-2}dt
$$
$$
{}
=
\int\limits_0^1
\frac{(1+t)^n-nt-1}{t^2}dt
=
\int\limits_0^1
\sum\limits_{s=1}^{n-1}
(n-s)(1+t)^{s-1}dt
=
\sum\limits_{s=1}^{n-1}
\frac{n-s}{s}(1+t)^s\Big|_0^1
=
\sum\limits_{s=1}^{n-1}
\frac{n-s}{s}(2^s-1)
$$
$$
{}=
\sum\limits_{s=1}^{n}
\frac{n-s}{s}2^s
-
\sum\limits_{s=1}^{n}
\frac{n-s}{s}
=
n\sum\limits_{s=1}^{n}
\frac{2^s}{s}
-
\sum\limits_{s=1}^{n}
2^s
-
n\sum\limits_{s=1}^{n}
\frac{1}{s}
+n
=
n\sum\limits_{s=1}^{n}
\frac{2^s}{s}
-
\frac{2(2^n-1)}{2-1}
-
n\sum\limits_{s=1}^{n}
\frac{1}{s}
+n.
$$
This equality and~(\ref{lemen_sum}) imply that
$$
\Upsilon(n)
\le
n\left(\frac{2^{n+1}}{n}
+\frac{2^{n+1}}{(n-3)^2}\right)
-2^{n+1}+2
-
n\sum\limits_{s=2}^{n}
\frac{1}{s}
=
\frac{2^{n+1}n}{(n-3)^2}
+2
-
n\sum\limits_{s=2}^{n}
\frac{1}{s}
$$
$$
{}
\le
\frac{2^{n+1}}{n}\frac{n^2}{(n-3)^2}
\le
\frac{2^{n+1}}{n}
\left(
1+\frac{15}{n-3}
\right).
$$
\hfill$\Box$

Now we can estimate $\bar{\mathcal{S}}_6$.

\begin{lemen}
\label{S_6}
$
\displaystyle{
\bar{\mathcal{S}}_6
\le
\frac{5\cdot2^{n+1}}{en(\ln n-1)}\left(1+\frac{15}{n-3}\right).
}
$
\end{lemen}

\textsl{Proof}.
Since $\bar{\mathcal{S}}_6=\sigma_1+\ldots+\sigma_4$, then
from Lemmas~\ref{sigma_12}, \ref{sigma_3}, \ref{sigma_4}
we obtain
$$
\bar{\mathcal{S}}_6
\le
\frac{C_0+C_0
}{e(\ln n-1)}
\sum\limits_{i_1=0}^{n-4}
\left(
\begin{array}{c}
n\\
n-i_1
\end{array}
\right)
\frac{1}{n-i_1-1}
+
\frac{3}{e(\ln n-1)}
\left(
\begin{array}{c}
n\\
n-3
\end{array}
\right)
\frac{1}{2}
+
\frac{2}{e(\ln n-1)}
\left(
\begin{array}{c}
n\\
n-2
\end{array}
\right)
$$
$$
{}
\le
\frac{5
}{e(\ln n-1)}
\sum\limits_{i_1=0}^{n-2}
\left(
\begin{array}{c}
n\\
n-i_1
\end{array}
\right)
\frac{1}{n-i_1-1}.
$$
This estimate and~(\ref{lemen_Ture_2}) imply the assertion of the lemma.
\hfill$\Box$

{\bf 8.3. The upper bounds for other sums.}
\begin{lemen}
\label{125}
The following inequalities hold:
$$
\bar{\mathcal{S}}_1+\bar{\mathcal{S}}_2+\bar{\mathcal{S}}_5
\le
\frac{2^{n+1}}{en(\ln n -\ln 2)}
\left(1+\frac{15}{n-3}\right)
\left(
1+\frac{e(\ln n-1)n(n+1)}{2^{n+1}}
\right),
$$
$$
\bar{\mathcal{S}}_3+\bar{\mathcal{S}}_4
\le\frac{2^{n+1}}{en(\ln n -\ln 2)}
\left(1+\frac{15}{n-3}\right)
\left(
1+\frac{e(\ln n-1)n^2}{2^{n+1}}
\right).
$$
\end{lemen}

\textsl{Proof}.
Recall that $r_3=0$. Then
$$
\bar{\mathcal{S}}_1+\bar{\mathcal{S}}_5+\bar{\mathcal{S}}_2
=\sum\limits_{i_2=0}^n
|l_i(\lambda)|\Big|_{i=(n-i_2,i_2,0)}.
$$
Consider $l_i(\lambda)$ for $i_2=0,1$:
\begin{equation}
\label{1}
|l_i(\lambda)|\Big|_{i=(n,0,0)}
=
\frac{\Gamma(n+\alpha_1)}{n!~\Gamma(\alpha_1)}
\le
\frac{\Gamma(n+1)}{n!~\Gamma(\alpha_1)}\le 1,
\end{equation}

\begin{equation}
\label{2}
|l_i(\lambda)|\Big|_{i=(n-1,1,0)}
=
\frac{\Gamma(n+\alpha_1)}{(n-1)!~\Gamma(1+\alpha_1)}
\frac{\Gamma(1+\alpha_2)}{\Gamma(\alpha_2)}
\le
\frac{\Gamma(n+1)}{(n-1)!~\Gamma(1+1)}
\alpha_2
\le n
\end{equation}
(we used Lemma~\ref{monot} to obtain the inequality~(\ref{2})).

Now we have to estimate the sum
$$
\sigma_0
\overset{def}{=}
\sum\limits_{i_2=2}^n
|l_i(\lambda)|\Big|_{i=(n-i_2,i_2,0)}
=
\sum\limits_{i_2=2}^n
\frac{\sin\pi\alpha_2}{\pi}
\frac{\Gamma(n+\alpha_1)}{(n-i_2)!~\Gamma(i_2+\alpha_1)}
\frac{\Gamma(1+\alpha_2)\Gamma(i_2-\alpha_2)}{i_2!}
$$
$$
{}=
\sum\limits_{i_2=2}^n
\frac{\sin\pi\alpha_2}{\pi}
\frac{\Gamma(n+1-\alpha_2-\alpha_3)}{(n-i_2)!~\Gamma(i_2+1-\alpha_2-\alpha_3)}
\frac{\Gamma(1+\alpha_2)\Gamma(i_2-\alpha_2)}{i_2!}.
$$
Lemma~\ref{monot} provides the estimate
$$
\sigma_0\le
\sum\limits_{i_2=2}^n
\frac{\sin\pi\alpha_2}{\pi}
\frac{\Gamma(n+1-\alpha_2)\Gamma(1+\alpha_2)}{(n+1)!}
\frac{(n+1)!~\Gamma(i_2-\alpha_2)}{i_2!~(n-i_2)!~\Gamma(i_2+1-\alpha_2)}.
$$
Using Lemmas~\ref{Ture_1}, \ref{ln},
\ref{Ture_2} we obtain
$$
\sigma_0\le
\sum\limits_{i_2=2}^n
\frac{\alpha_2(1-\alpha_2)2^{\alpha_2}}{(n+1)^{1+\alpha_2}}
\frac{(n+1)!~\Gamma(i_2-\alpha_2)}{i_2!~(n-i_2)!~\Gamma(i_2+1-\alpha_2)}
\le
\sum\limits_{i_2=2}^n
\frac{1-\alpha_2}{e(\ln n -\ln 2)}
\frac{n!}{i_2!~(n-i_2)!}~\frac{1}{i_2-\alpha_2}
$$
\begin{equation}
\label{3}
\le
\frac{1-\alpha_2}{e(\ln n -\ln 2)}
\sum\limits_{i_2=2}^n
\left(
\begin{array}{c}
n\\
i_2
\end{array}
\right)
\frac{1}{i_2-1}
{}\le
\frac{1}{e(\ln n -\ln 2)}
\frac{2^{n+1}}{n}\left(1+\frac{15}{n-3}\right).
\end{equation}
So, the first inequality of the lemma
follows from~(\ref{1}), (\ref{2}), (\ref{3}).
The second inequality is proved in a similar way.
\hfill$\Box$

{\bf 8.4. Completion of the proof of Theorem~\ref{teo_osn}}.
It follows from~(\ref{ddd}), (\ref{vvv}) and
Lemmas~\ref{S_6}, \ref{125} that
\begin{equation}
\label{t1t1}
\Lambda_n
\le
(7+\bar{\mu}_1+\bar{\mu}_2)
\frac{2^{n+1}}{en(\ln n -\ln 2)}
\left(1+\frac{15}{n-3}\right),
\end{equation}

where

\begin{equation}
\label{t2t2}
\bar{\mu}_1
\le\frac{e(\ln n-1)n(n+1)}{2^{n+1}}
+
\frac{e(\ln n-1)n^2}{2^{n+1}}
=
\frac{en(2n\ln n+\ln n-2n-1)}{2^{n+1}}
\le
\frac{en^2\ln n}{2^n},
\end{equation}

$$
\bar{\mu}_2
\le2^{2n/3}(10+2\ln n)
\left(
\frac{2^{n+1}}{en(\ln n -\ln 2)}
\left(1+\frac{15}{n-3}\right)
\right)^{-1}
$$
$$
{}
\le
\frac{en}{2^{n/3}}
(5+\ln n)(\ln n-\ln 2)
\frac{n-3}{n+12}
\le
\frac{en(\ln n)^2}{2^{n/3}}
\left(
\frac{5-\ln 2}{\ln n}+1-\frac{\ln 2}{(\ln n)^2}
\right).
$$
The point $t^*={\displaystyle \frac{10\ln 2}{5-\ln 2}}$ %\in[\ln 4,\ln 5]$
is the only maximum point of the function
$
\displaystyle{
\chi(t)=\frac{5-\ln 2}{t}+1-\frac{5\ln 2}{t^2}
}$, $t>0$,
and
$
\chi(t^*)<3.
$
Then
\begin{equation}
\label{t3t3}
\bar{\mu}_2
<
\frac{3en(\ln n)^2}{2^{n/3}}.
\end{equation}
The estimates~(\ref{t1t1}), (\ref{t2t2}), (\ref{t3t3}) imply~(\ref{teo_teo_osn}).
The proof of Theorem~\ref{teo_osn} is complete.

\vspace{5mm}
{\bf Acknowledgement.}
The author was funded by ISF grant 519/17.

\noindent
%Байдакова Наталия Васильевна\\ %\hfill Поступила  01.07.2020\\
%Институт математики и механики им. Н.Н.Красовского УрО РАН, г. Екатеринбург\\
%Уральский федеральный университет, г. Екатеринбург\\

%\newpage

~

%\vspace{6cm}

%$$
%\sum\limits_{s=1}^{n}
%\frac{2^s}{s}
%\le
%\frac{2^{n+1}}{n}+O\left(\frac{2^{n+1}}{n^2}\right),
%$$
%где $n\to \infty$. Например,
%$$
%\sum\limits_{s=1}^{n}
%\frac{2^s}{s}
%\le
%\frac{2^{n+1}}{n}+\frac{2^{n+1}}{(n-3)^2}
%$$

\end{document}